\documentclass[11pt,a4paper]{article}
\usepackage[T1]{fontenc}
\usepackage[a4paper, margin=1in]{geometry}
\usepackage{mathtools}
\usepackage{amssymb}
\usepackage{amsthm}
\usepackage{xcolor}
\usepackage[english]{babel}
\usepackage{hyperref}
\usepackage{amsmath}
\usepackage{mathrsfs}
\usepackage{cleveref}
\usepackage{titlesec}

\titleformat{\section}
[hang]
{\normalfont\Large\bfseries}
{\thesection.}
{0.4em}
{}

\titleformat{\subsection}
[hang]
{\normalfont\large\bfseries}
{\thesubsection.}
{0.3em}
{}

\definecolor{myred}{RGB}{255,0,0}

\crefname{equation}{equation}{equations}
\Crefname{equation}{Equation}{Equations}

\hypersetup{
	colorlinks=true,
	linkcolor=myred,
	citecolor=myred,
	urlcolor=myred
}

\allowdisplaybreaks

\theoremstyle{plain}

\newtheorem{thm}{Theorem}[section]
\newtheorem{lem}{Lemma}[section]
\newtheorem{prop}{Proposition}[section]

\newtheorem{defi}{Definition}[section]
\newtheorem{remark}{Remark}[section]

\Crefname{section}{Section}{Sections}
\Crefname{subsection}{Subsection}{Subsections}
\numberwithin{equation}{section}

\newcommand{\XT}[1]{\mathbb{L}^{\infty,p}_{2,r}(#1)}

\title{\textbf{Global well-posedness for an inhomogeneous nonlinear Schr\"odinger equation with multiplicative noise}}
\author{
	Sandip Roy\footnote{Department of Mathematics, Indian Institute of Technology Indore, Indore 453552, India}\textsuperscript{\,\,}\textsuperscript{\,,}\footnote{Corresponding author\newline
		E-mail addresses: phd2301141004@iiti.ac.in (S. Roy), 
		debopriya@iiti.ac.in (D. Mukherjee)} , 
	Debopriya Mukherjee\footnotemark[1]}

\begin{document}
\maketitle

\begin{abstract}
\noindent We establish the existence and uniqueness of a global mild solution in
$L^{q}\bigl(\Omega; C([0,T];L^{2}(\mathbb{R}^{d}))\cap L^{p}(0,T;L^{r,2}(\mathbb{R}^{d}))\bigr)$,
to the inhomogeneous stochastic nonlinear Schr\"odinger equation with a linear Stratonovich noise and a scalar mass-dependent nonlinear It\^o noise in the mass-subcritical regime $0<b<\min\{2,d\}$, $0<\sigma<\frac{2-b}{d}$, in every space dimension $d\ge1$. The singular coefficient $|x|^{-b}$ belongs to no Lebesgue space over $\mathbb{R}^{d}$. This leads us to work on the whole space in the Lorentz scale, using $|x|^{-b}\in L^{d/b,\infty}(\mathbb{R}^{d})$ together with the Lorentz refinement of the Strichartz estimates. This helps us to bypass the standard decomposition of $\mathbb{R}^d$ in a ball and its complement and leads to a single choice of admissible pair. The argument of the global solution follows from an exact mass identity. The truncated equation satisfies a deterministic ordinary differential equation for $\|u(t)\|_{L^{2}}^{2}$, which is then combined with a subcritical bootstrap to remove the truncation.
\end{abstract}


\section{Introduction}\label{secintro}

\noindent Consider the inhomogeneous stochastic nonlinear Schr\"odinger equation
\begin{equation}\label{maineq}
	\begin{cases*}
		du(t)=-\big[\mathrm{i} Au(t) +\mathrm{i}\,\mathcal{N}(u(t))+ \beta u(t)\big]dt
		- \mathrm{i} B(u(t)) \circ d\mathcal{W}_1(t)
		- \mathrm{i} G(u(t)) d\mathcal{W}_2(t), & $t>0$,\\
		u(0)= u_0,
	\end{cases*}
\end{equation}
where
\begin{align}\label{defN}
	\mathcal{N}(u):=\lambda\,|x|^{-b}\,|u|^{2\sigma}u,
	\qquad \lambda\in\mathbb{R},\quad b>0,\quad \sigma>0,
\end{align}
$A=-\Delta$, $\beta>0$ is a damping parameter, and $u_{0}\in L^{\infty}(\Omega;L^{2}(\mathbb{R}^{d};\mathbb{C}))$.
The operator $B$ generates a linear Stratonovich noise and $G$ a nonlinear It\^o noise; we refer \Cref{Ass} for more details.

The deterministic equation obtained by switching off both noises,
\begin{align}\label{INLSdeter}
	\mathrm{i}\partial_{t}u+\Delta u=\lambda|x|^{-b}|u|^{2\sigma}u,
\end{align}
is the \emph{inhomogeneous nonlinear Schr\"odinger} (INLS) equation. This models the
propagation of a laser beam through a medium whose density has been modulated on purpose,
the factor $|x|^{-b}$ encoding the resulting spatial inhomogeneity of the nonlinear response. \Cref{INLSdeter} it also arises as the asymptotic profile of a polynomially decaying potential. In such media the modulation is never exactly reproducible, and thermal and mechanical fluctuations of the
channel act on the beam as a random perturbation of the refractive index. The precise reason behind the choice of the two noises in \eqref{maineq} modelis: (i) the Stratonovich term $-\mathrm{i}B(u)\circ d\mathcal{W}_{1}$ is the conservative (phase) noise generated by fluctuations of the linear refractive index, and (ii) $-\mathrm{i}G(u)d\mathcal{W}_{2}$ is a scalar noise whose amplitude depends on the global mass rather than on the pointwise intensity. The damping $\beta u$ accounts for the loss of the fibre. The equation is called defocusing INLS if $\lambda=1$ and focusing if $\lambda=1$. From this point of view the stochastic problem is not a formal generalisation of \eqref{INLSdeter} but the physically natural formulation. Since the latter noise is written in It\^o form, it is phase-orthogonal at the martingale level but contributes positively to the mass through its quadratic variation. This convention is accounted in the model.

\subsection{Literature review}
As best of our knowledge, there is no work on stochastic INLS in the literature. Though, in the deterministic set up (equation \eqref{INLSdeter}), it has been investigated by several authors over the last few decades. For the case $b=0$, equation \eqref{maineq} reduces to the classical nonlinear Schr\"odinger equation, which arises in many areas with a wide range of physical applications, for eg: quantum mechanics, nonlinear optics, plasma physics, Bose-Einstein condensates, hydrodynamics, molecular biology, see \cite{ agrawal2000nonlinear, PhysRevE.75.066607, PhysRevA.60.1507, chen1984introduction, karjanto2019nonlinear, shankar2012principles, sulem2007nonlinear, zhang2014stochastic}. The classical nonlinear Schr\"odinger equation has been widely studied by several authors over the last few decades in both deterministic, see \cite{T_Kato_1987, MR1162728, MR1055532, MR460850, MR1691575, cazenave} and stochastic settings, see \cite{MR3705702, MR3215081, MR3474409, Hornung_2020, MR3232027, MR3859442, MR4249116, MR4678996, deBourd+Debussche_2003}, in terms of the well-posedness and the long-term behavior of solutions, including stability and other qualitative properties. 

A more general form of \eqref{INLSdeter}, namely
\begin{align}\label{GeneralINLS}
	\mathrm{i}\partial_{t}u+\Delta u+V(x)|u|^{2\sigma}u=0,
\end{align}
where $V(x)$ is a potential, has been extinsively studied by several authors, see \cite{Merle_1996, Berge_2000, GFibich_WXP_2003, Debourd_RFukuizumi_2005, YLiu_XPW_KWang_2006, RPierre_SJeremie_2011}. For other interesting physical applications of \eqref{GeneralINLS}, one can see \cite{gill2000optical, liu1994laser, borovkova2011algebraic, borovkova2012solitons} and the references therein. In perticular, with $V(x)=|x|^{-b},\,\, b>0$, that is equation \eqref{INLSdeter} has been studied intensively. Genoud and Stuart \cite{GenoudStuart2008}
established local well-posedness in $H^{1}(\mathbb{R}^{d})$ for $0<b<\min\{2,d\}$ by an
argument that does not use Strichartz estimates. Guzm\'an \cite{Guzman2017} obtained local well-posedness
in $H^{s}$ via the contraction mapping principle based on Strichartz estimates, with parameter
	restrictions depending on $s$ and the dimension. The advantage is that the solution then lies in
$L^{p}_{\mathrm{loc}}(\mathbb{R};W^{1,r})$
for every admissible pair, which is required for scattering. In $L^{2}(\mathbb{R}^{d})$,
An and Kim \cite{AnKim2021} proved local well-posedness (and small-data global well-posedness)
in the whole range $0<b<\min\{2,d\}$, $0<2\sigma\le\frac{4-2b}{d}$, including the mass-critical
case, by working with Strichartz estimates in Lorentz spaces. The mass-critical exponent is
$\sigma=\frac{2-b}{d}$. For more references, we cite \cite{Genoud_2012, An_Kim2021, An_Kim2023, JChen_2010, Farah_2016, Kim_Lee_Seo_2021, Farah_Guzman_2017, Dinh_2018}. 

\subsection{Novelties and difficulties}\label{subsecnovelties}

To the best of our knowledge, the stochastic INLS equation \eqref{maineq}, with its
singular spatial coefficient and the specified Stratonovich and nonlinear It\^o noises,
has not previously been studied in the $L^2$ setting. Its analysis presents three
linked difficulties.

\smallskip
\noindent\textbf{(1) Singular weight and a single Strichartz pair.}
Usual H\"older inequalities cannot handle the weight on all of $\mathbb R^d$, since
\begin{align}\label{nolebesgue}
|x|^{-b}\notin L^{\rho}(\mathbb R^d)\quad\text{for every }1\le\rho\le\infty.
\end{align}
Instead of decomposing space into a neighbourhood of the origin and its complement, we use $|x|^{-b}\in L^{d/b,\infty}(\mathbb R^d)$ and the Lorentz-H\"older inequality, as in the deterministic approach of \cite{AnKim2021}. Requiring the nonlinear estimate in the dual of the same spatial Strichartz space,
$1/r'=b/d+(2\sigma+1)/r$, together with Schr\"odinger admissibility, determines the pair
\begin{align}\label{intropr}
r=\frac{2(\sigma+1)d}{d-b},\qquad
p=\frac{4(\sigma+1)}{\sigma d+b}.
\end{align}
For this pair (\Cref{lemexponents}), the spatial and space-time estimates are
\begin{align}\label{introest}
\|\mathcal N(u)\|_{L^{r',2}}\le C\|u\|_{L^{r,2}}^{2\sigma+1},
\end{align}
\begin{align}\label{introtheta}
\|\mathcal N(u)\|_{L^{p'}(0,T;L^{r',2})}
\le CT^{\theta}\|u\|_{L^p(0,T;L^{r,2})}^{2\sigma+1},
\quad\theta=\frac{2-b-\sigma d}{2}>0.
\end{align}
The positive time gain is exactly the mass-subcritical condition
$\sigma<(2-b)/d$ (\Cref{propspacetime}).

\smallskip
\noindent\textbf{(2) Lorentz-valued stochastic convolution.}
The fixed point is constructed in $C([0,T];H)\cap L^p(0,T;L^{r,2}(\mathbb R^d))$, whereas standard stochastic Strichartz estimates are formulated in Lebesgue spaces \cite{MR3232027,Hornung_2020,MR4249116}. We show that the relevant Lorentz-valued space admits an equivalent $2$-uniformly smooth norm
(\Cref{lem2smooth}). Also, Banach-space stochastic integration then yields a stochastic Strichartz estimate, including constants independent of $T$ (\Cref{stochStrichartz}).

\smallskip
\noindent\textbf{(3) Globalization at the $L^2$ level.} The nonlinear term requires higher moments, making it difficult to obtain estimates directly in $L^q(\Omega)$. A cut-off depending on the full
Strichartz norm gives a globally Lipschitz truncated problem. For its solutions, a Yosida regularisation justifies the mass identity, i.e., the nonlinear contribution vanishes by the $\lambda$ and the weight,
\begin{align}\label{introcancel}
(-\mathrm{i}\mathcal N(u),u)_H
=\operatorname{Re}\!\left(-\mathrm{i}\lambda
\int_{\mathbb R^d}|x|^{-b}|u|^{2\sigma+2}\,dx\right)=0.
\end{align}
The linear Stratonovich correction cancels its quadratic-variation
term, and both phase-noise martingales vanish. The resulting
mass identity gives a pathwise bound independent of the cut-off.
A localized subcritical bootstrap then removes the truncation and establishes global well-posedness (\Cref{massidentity,mainthm}).

\subsection{Organization of the paper}
\noindent The paper is organised as follows. In \Cref{secprelim}, we introduce the notation, state the assumptions and the main result, and converts \eqref{maineq} into complete It\^o form. In \Cref{secnonlinear}, we establish the new nonlinear estimates \eqref{introest} and \eqref{introtheta} together with their corresponding Lipschitz estimates. This section form the analytical foundation of the paper. In \Cref{sectruncated}, we introduce the truncated equation and global well-posedness is established in \Cref{secwellposedtrunc} by a fixed point argument followed by an induction. Finally, in \Cref{secremoval}, we establish the mass identity, remove the truncation and complete the proof of \Cref{mainthm}.

\section{Preliminaries and the statement of the main result}\label{secprelim}
In this section, we introduce some basic notations and recall some preliminary results. Next, we state the main result of this article.
\subsection{Notation}\label{subsecnotation}
We consider a filtered probability space $\mathfrak{U}=(\Omega, \mathcal{F}, \mathbb{F}, \mathbb{P})$,
where $\mathbb{F} = \{ \mathcal{F}_t \}_{t \geq 0}$, satisfying the usual conditions:
\begin{enumerate}
	\item $\mathbb{P}$ is complete on $(\Omega, \mathcal{F})$,
	\item for each $t \ge 0$, $\mathcal{F}_t$ contains all $(\mathcal{F}, \mathbb{P})$-null sets,
	\item the filtration $\mathbb{F}$ is right continuous, that is,
	$\mathcal{F}_t = \bigcap_{s > t} \mathcal{F}_s$ for all $t \ge 0$.
\end{enumerate}
We denote by $L^p(\mathbb{R}^d;\mathbb{C})$, $1\le p<\infty$, the Lebesgue space of (equivalence classes of)
complex-valued measurable functions whose $p$-th power is integrable, and we write
\begin{align*}
	H:=L^2(\mathbb{R}^d;\mathbb{C}),
\end{align*}
equipped with the real inner product
\begin{align*}
	(u,v)_H=\operatorname{Re}\int_{\mathbb{R}^d}u(x)\overline{v(x)}\,dx
\end{align*}
and induced norm $\|\cdot\|_{H}$. We write $\langle u,v\rangle:=\int_{\mathbb{R}^{d}}u\overline{v}\,dx$
for the Hermitian inner product, so that $(u,v)_{H}=\operatorname{Re}\langle u,v\rangle$.
For $1\le p,\varrho\le\infty$ and $T>0$, $L^p(0,T;X)$ denotes the usual Bochner space of strongly
measurable $X$-valued functions on $[0,T]$, and $C([0,T];H)$ the Banach space of $H$-valued continuous
functions with the supremum norm. For a Banach space $X$ and $1\le q<\infty$, $L^q(\Omega;X)$ is the
space of strongly measurable $X$-valued random variables $f$ with $\mathbb{E}\|f\|_X^q<\infty$, normed by
$\|f\|_{L^{q}(\Omega;X)}=(\mathbb{E}\|f\|^{q}_{X})^{1/q}$. We write $\mathcal{L}_{2}(Y,H)$ for the space
of Hilbert-Schmidt operators from $Y$ to $H$ and $\mathscr{L}(H)$ for the bounded linear operators on $H$.
Finally, $\{S(t)\}_{t\in\mathbb{R}}$ denotes the unitary group generated by $-\mathrm{i}A=\mathrm{i}\Delta$, i.e. $S(t)=e^{\mathrm{i}t\Delta}$, which is the group appearing in the mild formulation of \eqref{maineq}.

\textbf{Lorentz spaces:}
For a measurable $f:\mathbb{R}^{d}\to\mathbb{C}$, let
$d_{f}(\alpha)=|\{x:|f(x)|>\alpha\}|$ be its distribution function and
$f^{*}(t)=\inf\{\alpha>0: d_{f}(\alpha)\le t\}$ its decreasing rearrangement. For
$0<\varrho<\infty$ and $0<s\le\infty$ the Lorentz space $L^{\varrho,s}(\mathbb{R}^{d})$ consists of
those $f$ for which
\begin{align}\label{lorentznorm}
	\|f\|_{L^{\varrho,s}}:=
	\begin{cases}
		\displaystyle\Bigl(\int_{0}^{\infty}\bigl(t^{1/\varrho}f^{*}(t)\bigr)^{s}\frac{dt}{t}\Bigr)^{1/s},
		& 0<s<\infty,\\[2ex]
		\displaystyle\sup_{t>0}\;t^{1/\varrho}f^{*}(t), & s=\infty,
	\end{cases}
\end{align}
is finite. For a subinterval $I\subseteq[0,T]$, we define the Banach space
\begin{align}\label{defXS}
	\mathcal{X}({I}):=C(\overline I;H)\cap L^{p}(I;L^{r,2}(\mathbb{R}^{d})),\qquad \mathcal{X}({t}):=\mathcal{X}({[0,t]}) .
\end{align}
with the norm
\begin{align}\label{normassum}
	\|u\|_{\mathcal{X}({I})}:=\|u\|_{L^{\infty}(I;H)}+\|u\|_{L^{p}(I;L^{r,2}(\mathbb{R}^{d}))}.
\end{align}

\subsection{Assumptions}\label{Ass}
\begin{itemize}
	\item[\textbf{(a)}] \textbf{The operator $A$.} $A=-\Delta$, the negative Laplacian, as an operator on $H$.
	\item[\textbf{(b)}] \textbf{The nonlinearity.} $\mathcal{N}(u)=\lambda|x|^{-b}|u|^{2\sigma}u$ with $\lambda\in\mathbb{R}$,
	\begin{align}\label{brange}
		0<b<\min\{2,d\},\qquad 0<\sigma<\frac{2-b}{d}.
	\end{align}
	The second condition is mass-subcriticality; the mass-critical exponent is $\sigma=\frac{2-b}{d}$.
	That $\lambda$ is \emph{real} is essential for the mass identity \eqref{introcancel}; the sign of
	$\lambda$ (focusing or defocusing) is irrelevant in the subcritical regime treated here.
	\item[\textbf{(c)}] \textbf{The noise and its coefficients.}
	\begin{itemize}
		\item[(i)] Let $Y_{1}$ be a separable real Hilbert space with orthonormal basis
		$\{e^{1}_{m}\}_{m\in\mathbb{N}}$ and let $\mathcal{W}_{1}$ be a cylindrical
		$\mathbb{F}$-Wiener process on $Y_{1}$,
		\begin{align}\label{BasisW}
			\mathcal{W}_1(t)=\sum_{m=1}^{\infty}e^{1}_{m}\beta^{1}_{m}(t),\qquad t\ge0,
		\end{align}
		with $\{\beta^{1}_{m}\}_{m\in\mathbb{N}}$ a sequence of independent standard real Wiener processes.
		The saturable noise is taken \emph{scalar}: we set $Y_{2}:=\mathbb{R}$ and write
		$\mathcal{W}_{2}=:\beta^{2}$ for a standard real $\mathbb{F}$-Wiener process independent of
		$\mathcal{W}_{1}$. Throughout we use the canonical isometric identification
		\begin{align}\label{L2R}
			\mathcal{L}_{2}(\mathbb{R},H)\;\cong\;H,\qquad \Phi\longmapsto\Phi(1),
		\end{align}
		under which $\|\Phi\|_{\mathcal{L}_{2}(\mathbb{R},H)}=\|\Phi(1)\|_{H}$.
		\item[(ii)] Let $B:H\to\mathcal{L}_{2}(Y_{1},H)$ be linear and set $B_{m}u:=B(u)e^{1}_{m}$, so that
		$B_{m}\in\mathscr{L}(H)$. We assume each $B_{m}$ is $\mathbb{C}$-linear and self-adjoint with respect to
		the \emph{Hermitian} inner product, i.e.
		\begin{align}\label{Bmselfadjoint}
			\langle B_{m}u,v\rangle=\langle u,B_{m}v\rangle,\qquad u,v\in H,
		\end{align}
		and that, in order for the Stratonovich correction to make sense,
		\begin{align}\label{BmH}
			\sum_{m=1}^{\infty}\|B_{m}\|^{2}_{\mathscr{L}(H)}<\infty.
		\end{align}
		\item[(iii)] By \eqref{L2R} the saturable coefficient is a map $G:H\to H$, and we take it of
		\emph{saturation type},
		\begin{align}\label{perticularg}
			G(u)=\tilde G\bigl(\|u\|^{2}_{H}\bigr)\,u,
		\end{align}
		where $\tilde G:[0,\infty)\to\mathbb{R}$ is measurable and bounded. Then, by \eqref{L2R},
		\begin{align}\label{LgG}
			\|G(u)\|_{\mathcal{L}_{2}(Y_{2},H)}
			=\bigl|\tilde G\bigl(\|u\|^{2}_{H}\bigr)\bigr|\,\|u\|_{H}
			\le\|\tilde G\|_{\infty}\|u\|_{H},\qquad u\in H,
		\end{align}
		so that $G$ is of linear growth with vanishing constant term. Under \eqref{L2R} the norm
		$\|\cdot\|_{\mathcal{L}_{2}(Y_{2},H)}$ is simply $\|\cdot\|_{H}$, we keep the former notation for
		uniformity with the $\mathcal{W}_{1}$-terms. We assume in addition that $G$ is Lipschitz,
		\begin{align}\label{LipG}
			\|G(u_{1})-G(u_{2})\|_{\mathcal{L}_{2}(Y_{2},H)}\le L_{G}\|u_{1}-u_{2}\|_{H},
			\qquad u_{1},u_{2}\in H,
		\end{align}
	\end{itemize}
\end{itemize}

\begin{remark}\label{remsatnl}
	Condition \eqref{LipG} is true when $\tilde G\in C^{1}([0,\infty))$ satisfies
	\begin{align}\label{Gsuff}
		\|\tilde G\|_{\infty}+2\sup_{y\ge0}\,y\,\bigl|\tilde G'(y)\bigr|<\infty ,
	\end{align}
	with $L_{G}$ bounded by the left-hand side of \eqref{Gsuff}. Indeed
	$$D\bigl[\tilde G(\|u\|^{2}_{H})u\bigr][h] =\tilde G(\|u\|^{2}_{H})h+2\tilde G'(\|u\|^{2}_{H})(u,h)_{H}\,u,$$
	 and
	$$\|2\tilde G'(\|u\|^{2}_{H})(u,h)_{H}u\|_{H}\le2\|u\|^{2}_{H}|\tilde G'(\|u\|^{2}_{H})|\,\|h\|_{H}.$$
	The class of saturation nonlinearities satisfying \eqref{perticularg}-\eqref{LipG} and \eqref{Gsuff}
	is large. Here, one may take $\tilde G(y)=\frac{y}{1+ay}$ ($a>0$), which describes photorefractive media \cite{PhysRevLett.15.1005}, or $\tilde G(y)=1-(1+y)^{-1/2}$, which describes narrow-gap semiconductors \cite{PhysRevB.61.10201}, or
	$\tilde G(y)=e^{-y}$, $\tilde G(y)=\frac{y(ay+2)}{(1+ay)^{2}}$, $\tilde G(y)=\frac{\log(1+ay)}{1+\log(1+ay)}$. In each case $y\tilde G'(y)$ is bounded on $[0,\infty)$.
\end{remark}

\begin{defi}[Admissible pair]\label{Admissible}
	A pair $(p,\varrho)\in[2,\infty]\times[2,\infty]$ is \emph{admissible} if
	\begin{align*}
		\frac{2}{p}=\frac{d}{2}-\frac{d}{\varrho},\qquad (p,\varrho,d)\neq(2,\infty,2).
	\end{align*}
\end{defi}

\subsection{Choice of exponents}
The choice of the spatial exponent $r$ is due to the singular weight
$|x|^{-b}\in L^{d/b,\infty}(\mathbb R^d)$ along with the dual Lorentz exponent, 
$1/r'=b/d+(2\sigma+1)/r$. Combining this identity with the
Schr\"odinger admissibility relation $2/p=d/2-d/r$ determines the pair
$(p,r)$ used below. Thus the choice is dictated by the nonlinear
Lorentz-H\"older estimate (\Cref{lemoneil}) and the linear Strichartz
estimates.
\begin{lem}[Admissibility]\label{lemexponents}
	Define
	\begin{align}\label{defpr}
		r:=\frac{2(\sigma+1)d}{d-b},\qquad p:=\frac{4(\sigma+1)}{\sigma d+b},
		\qquad \theta:=1-\frac{2(\sigma+1)}{p} .
	\end{align}
	Then
	\begin{enumerate}
		\item[(i)] $2<r<\infty$; moreover $r<\frac{2d}{d-2}$ when $d\ge3$, so that $(p,r)$ is an admissible pair
		in the sense of \Cref{Admissible};
		\item[(ii)] $\dfrac{2}{p}=\dfrac{d}{2}-\dfrac{d}{r}$;
		\item[(iii)] $p>2(\sigma+1)>2$ and
		\begin{align}\label{deftheta}
			\theta=\frac{2-b-\sigma d}{2}\in(0,1).
		\end{align}
	\end{enumerate}
	Conversely, $\theta>0$ holds if and only if $\sigma<\frac{2-b}{d}$; that is, the time gain in
	\Cref{propspacetime} is available precisely in the mass-subcritical regime.
\end{lem}

\begin{proof}
	(i) Since $b<d$ we have $\frac{d}{d-b}>1$, whence $r>2(\sigma+1)>2$; and $r<\infty$. For $d\ge3$,
	\begin{align*}
		r<\frac{2d}{d-2}
		\iff (\sigma+1)(d-2)<d-b
		\iff \sigma(d-2)<2-b
		\iff \sigma<\frac{2-b}{d-2},
	\end{align*}
	and the last inequality holds because $\sigma<\frac{2-b}{d}<\frac{2-b}{d-2}$. For $d=2$ the requirement
	is $2\le r<\infty$, and for $d=1$ it is $2\le r\le\infty$ both hold, since $2<r<\infty$ was already
	shown. Since $r<\infty$, the excluded endpoint $(2,\infty,2)$ does not occur.
	
	(ii) We compute
	\begin{align*}
		\frac{d}{2}-\frac{d}{r}=\frac{d}{2}-\frac{d-b}{2(\sigma+1)}
		=\frac{d(\sigma+1)-(d-b)}{2(\sigma+1)}=\frac{\sigma d+b}{2(\sigma+1)}
		=\frac{2}{p},
	\end{align*}
	by the definition of $p$.
	
	(iii) From \eqref{defpr},
	$\frac{2(\sigma+1)}{p}=\frac{2(\sigma+1)(\sigma d+b)}{4(\sigma+1)}=\frac{\sigma d+b}{2}$,
	so $\theta=1-\frac{\sigma d+b}{2}=\frac{2-b-\sigma d}{2}$. Since $b,\sigma d>0$ we get $\theta<1$,
	and $\theta>0$ is equivalent to $\sigma d+b<2$, i.e.\ $\sigma<\frac{2-b}{d}$, which is \eqref{brange}.
	Finally $\theta>0$ is the same as $p>2(\sigma+1)$.
\end{proof}

\noindent The following identity, which is a consequence of the definition of $r$, is used repeatedly:
\begin{align}\label{keyidentity}
	1-\frac{b}{d}=\frac{d-b}{d}=\frac{2(\sigma+1)}{r},
	\qquad\text{equivalently}\qquad
	\frac{1}{r'}-\frac{b}{d}=\frac{2\sigma+1}{r}.
\end{align}
Indeed, $\frac{1}{r'}-\frac{b}{d}=1-\frac1r-\frac bd=\frac{2(\sigma+1)}{r}-\frac1r=\frac{2\sigma+1}{r}$.
\subsection{The statement of the main result}
In this subsection, we state the main result of the paper. Before stating the main result, we define the notion of the solution in the mild sense.
\begin{defi}[Mild solution]\label{milditodefi}
	Fix $q\ge2$ and let $(p,r)$ be an admissible pair defined in \eqref{defpr}. Given $u_{0}\in L^{\infty}(\Omega,\mathcal{F}_{0};H)$,
	a \emph{mild solution} of \eqref{itoform} is an $\{\mathcal{F}_{t}\}_{t\ge0}$-adapted
	and progressively measurable process
	$u\in L^{q}\bigl(\Omega;\mathcal{X}({T})\bigr)$ such that
	\begin{align}\label{milditoeq}
		u(t)&=S(t)u_{0}-\int_{0}^{t}S(t-s)\bigl[\mathrm{i}\mathcal{N}(u(s))+\beta u(s)+\mathfrak{b}(u(s))\bigr]ds\nonumber\\
		&\quad-\mathrm{i}\int_{0}^{t}S(t-s)\bigl[B(u(s))d\mathcal{W}_1(s)+G(u(s))d\mathcal{W}_2(s)\bigr]
	\end{align}
	for every $t\in[0,T]$ $\mathbb P$-a.s.
\end{defi}

\begin{thm}\label{mainthm}
	Let $d\ge1$, $0<b<\min\{2,d\}$ and $0<\sigma<\frac{2-b}{d}$, and let $(p,r)$ be given by
	\eqref{intropr}. Let $q\ge2$, $T>0$ and $u_{0}\in L^{\infty}(\Omega,\mathcal{F}_{0};H)$.
	Under the assumptions of \Cref{Ass}, there exists a unique global mild solution
	\begin{align*}
		u\in L^{q}\bigl(\Omega; C([0,T];H)\cap L^{p}(0,T;L^{r,2}(\mathbb{R}^{d}))\bigr)
	\end{align*}
	of \eqref{maineq}, in the sense of \Cref{milditodefi}. Moreover the mass obeys the pathwise bound
	\begin{align}\label{mainmass}
		\sup_{t\in[0,T]}\|u(t)\|^{2}_{H}
		\le \|u_{0}\|^{2}_{H}\,e^{\bigl(\|\tilde G\|^{2}_{\infty}-2\beta\bigr)_{+}T},
		\qquad \mathbb{P}\text{-a.s.}
	\end{align}
	In particular, if the damping dominates the saturable noise in the sense that
	$2\beta\ge\|\tilde G\|^{2}_{\infty}$, then $t\mapsto\|u(t)\|_{H}$ is non-increasing and
	$\|u(t)\|_{H}\le\|u_{0}\|_{H}$ for every $t\in[0,T]$, $\mathbb{P}$-a.s.
\end{thm}

\subsection{Conversion to complete It\^o form}\label{subsecito}
The Stratonovich noise in \eqref{maineq} is convenient for the physics but not for the mathematical estimates, so we convert it into the It\^o form. Using \eqref{BasisW}, we get
\begin{align*}
	-\mathrm{i}B(u(t))\circ d\mathcal{W}_1(t)
	=\sum_{m\in\mathbb{N}}-\mathrm{i}B(u(t))e^{1}_{m}\circ d\beta^{1}_{m}(t)
	=:\sum_{m\in\mathbb{N}}\varsigma_{m}(u)\circ d\beta^{1}_{m}(t),
\end{align*}
where $\varsigma_{m}(u)=-\mathrm{i}B_{m}u$. Its Fr\'echet derivative in the direction $h$ is
$D\varsigma_{m}(u)[h]=-\mathrm{i}B_{m}h$, so
\begin{align*}
	-\mathrm{i}B(u(t))\circ d\mathcal{W}_1(t)
	&=-\mathrm{i}B(u(t))\,d\mathcal{W}_1(t)+\frac{1}{2}\sum_{m=1}^{\infty}D\varsigma_{m}(u)[\varsigma_{m}(u)]\,dt\\
	&=-\mathrm{i}B(u(t))\,d\mathcal{W}_1(t)-\frac{1}{2}\sum_{m=1}^{\infty}B_{m}^{2}u(t)\,dt .
\end{align*}
Hence \eqref{maineq} takes the complete It\^o form
\begin{equation}\label{itoform}
	du(t)=-\big[\mathrm{i}Au(t)+\mathrm{i}\mathcal{N}(u(t))+\beta u(t)+\mathfrak{b}(u(t))\big]dt
	-\mathrm{i}B(u(t))d\mathcal{W}_1(t)-\mathrm{i}G(u(t))d\mathcal{W}_2(t),\quad t\ge0,
\end{equation}
with the It\^o-Stratonovich correction
\begin{align}\label{defb}
	\mathfrak{b}(u):=\frac{1}{2}\sum_{m=1}^{\infty}B_{m}^{2}u .
\end{align}
By \eqref{BmH}, $\mathfrak{b}\in\mathscr{L}(H)$ with $\|\mathfrak{b}\|_{\mathscr{L}(H)}\le\frac12\sum_{m}\|B_{m}\|^{2}_{\mathscr{L}(H)}$.

\section{ Estimates for the nonlinear term}\label{secnonlinear}
This section develops the key analytical framework underlying the main result of the paper. In particular, we establish several crucial estimates for the nonlinear term that will be repeatedly used in the subsequent analysis. Throughout, we assume that $d\ge1$ and $b,\sigma$ satisfy \eqref{brange}. The constant $C$ denotes a generic positive constant, may change from line to line and depend only on $d,b,\sigma$ unless otherwise indicated.

\subsection{The pointwise-in-time estimate}
The weight $|x|^{-b}$ is not in any global Lebesgue space $L^\rho(\mathbb R^d)$, so ordinary H\"older estimates do not directly control $\mathcal N(u)$ using a single Lebesgue exponent. In contrast, $|x|^{-b}\in L^{d/b,\infty}(\mathbb R^d)$, and the exponent identity \eqref{keyidentity} allows Lorentz-H\"older to control the product in $L^{r',2}$ when $u\in L^{r,2}$. We first establish this spatial bound and then its Lipschitz counterpart, together they provide the pointwise estimates needed for the space-time bounds in
\Cref{propspacetime}.
\begin{lem}\label{lemnonlinearspace}
	Let $r$ be as in \eqref{defpr}. There is $C=C(d,b,\sigma,\lambda)>0$ such that
	\begin{align}\label{nonlinearspacebound}
		\bigl\|\mathcal{N}(u)\bigr\|_{L^{r',2}(\mathbb{R}^{d})}\le C\,\|u\|^{2\sigma+1}_{L^{r,2}(\mathbb{R}^{d})}
	\end{align}
	for every $u\in L^{r,2}(\mathbb{R}^{d})$.
\end{lem}

\begin{proof}
	We choose $q:=\frac{r}{2\sigma+1}$. Then, by \eqref{keyidentity}, we have
	\begin{align}\label{holderexponents}
		\frac{1}{r'}=\frac{b}{d}+\frac{1}{q}.
	\end{align}
	We check the condition for \Cref{lemoneil} with $\varrho_{1}=\frac{d}{b}$, $\varrho_{2}=q$,
	$s_{1}=\infty$, $s_{2}=2$, $\varrho=r'$, $s=2$. Since $b<d$, we have $\frac{d}{b}>1$. Also, $r>2(\sigma+1)>2\sigma+1$, so $q=\frac{r}{2\sigma+1}>1$. The condition $\frac{1}{r'}<1$ follows from $r<\infty$.

	Therefore, using \Cref{lemoneil}, \Cref{lemweight} and then \Cref{lemlorentzbasic}(iii) with $\alpha=2\sigma+1$,
	\begin{align*}
		\bigl\|\mathcal{N}(u)\bigr\|_{L^{r',2}}
		&=|\lambda|\,\bigl\|\,|x|^{-b}\,|u|^{2\sigma}u\,\bigr\|_{L^{r',2}}\\
		&\le C\,\bigl\||x|^{-b}\bigr\|_{L^{\frac{d}{b},\infty}}\,
		\bigl\||u|^{2\sigma}u\bigr\|_{L^{q,2}}
		= C\,\omega_{d}^{b/d}\,\|u\|^{2\sigma+1}_{L^{(2\sigma+1)q,\,2(2\sigma+1)}} .
	\end{align*}
	By the definition of $q$ we have $(2\sigma+1)q=r$, and since $2\le 2(2\sigma+1)$,
	\Cref{lemlorentzbasic}(ii) gives
	$\|u\|_{L^{r,2(2\sigma+1)}}\le C\|u\|_{L^{r,2}}$. This proves \eqref{nonlinearspacebound}.
\end{proof}

\begin{lem}\label{lemlipspace}
	There is $C=C(d,b,\sigma,\lambda)>0$ such that, for all $u,v\in L^{r,2}(\mathbb{R}^{d})$,
	\begin{align}\label{lipofnonlinearinspace}
		\bigl\|\mathcal{N}(u)-\mathcal{N}(v)\bigr\|_{L^{r',2}(\mathbb{R}^{d})}
		\le C\Bigl(\|u\|^{2\sigma}_{L^{r,2}}+\|v\|^{2\sigma}_{L^{r,2}}\Bigr)\|u-v\|_{L^{r,2}} .
	\end{align}
\end{lem}

\begin{proof}
	We first establish the following pointwise inequality
	\begin{align}\label{pointwise}
		\bigl||z_{1}|^{2\sigma}z_{1}-|z_{2}|^{2\sigma}z_{2}\bigr|
		\le C(\sigma)\bigl(|z_{1}|^{2\sigma}+|z_{2}|^{2\sigma}\bigr)|z_{1}-z_{2}|,
		\qquad z_{1},z_{2}\in\mathbb{C}.
	\end{align}
	This estimate follows from the mean value theorem applied to $z\mapsto|z|^{2\sigma}z$. Hence, pointwise in $x$, we have
	\begin{align*}
		\bigl|\mathcal{N}(u)-\mathcal{N}(v)\bigr|\le C|x|^{-b}\bigl(|u|^{2\sigma}+|v|^{2\sigma}\bigr)|u-v| .
	\end{align*}
	We now apply \Cref{lemoneil} twice. Choose
	\begin{align*}
		\frac{1}{m}:=\frac{b}{d}+\frac{2\sigma}{r},\qquad
		\frac{1}{\tilde s}:=\frac{2\sigma}{2(2\sigma+1)}=\frac{\sigma}{2\sigma+1}.
	\end{align*}
	Using \eqref{keyidentity}, we get $$\frac{1}{m}<\frac{b}{d}+\frac{2\sigma+1}{r}=\frac{1}{r'}<1,$$
  consequently $m>1$. Since $r>2(\sigma+1)>2\sigma$, we have $\frac{r}{2\sigma}>1$. Thus the pairs
	$(\varrho_{1},s_{1})=(\frac{d}{b},\infty)$ and
	$(\varrho_{2},s_{2})=\bigl(\frac{r}{2\sigma},\frac{2(2\sigma+1)}{2\sigma}\bigr)$ satisfy the condition of \Cref{lemoneil}. Using \Cref{lemoneil} together with
	\Cref{lemweight} and \Cref{lemlorentzbasic}(ii),(iii), we infer
	\begin{align}\label{mstilenorm}
	\bigl\||x|^{-b}\bigl(|u|^{2\sigma}+|v|^{2\sigma}\bigr)\bigr\|_{L^{m,\tilde s}}\le C\omega_{d}^{b/d} \Bigl(\|u\|^{2\sigma}_{L^{r,2(2\sigma+1)}}+\|v\|^{2\sigma}_{L^{r,2(2\sigma+1)}}\Bigr)
	\le C\Bigl(\|u\|^{2\sigma}_{L^{r,2}}+\|v\|^{2\sigma}_{L^{r,2}}\Bigr).
	\end{align}
 The condition of \Cref{lemoneil} with
	$(\varrho_{1},s_{1})=(m,\tilde s)$ and $(\varrho_{2},s_{2})=(r,2(2\sigma+1))$, is justified since
	$1<m<\infty$, $1<r<\infty$,
	\begin{align*}
		\frac{1}{m}+\frac{1}{r}=\frac{b}{d}+\frac{2\sigma+1}{r}=\frac{1}{r'}<1,
		\qquad
		\frac{1}{\tilde s}+\frac{1}{2(2\sigma+1)}=\frac{2\sigma+1}{2(2\sigma+1)}=\frac{1}{2}.
	\end{align*}
A second application of \Cref{lemoneil} yields
	\begin{align*}
		\bigl\|\mathcal{N}(u)-\mathcal{N}(v)\bigr\|_{L^{r',2}}
		\le C\bigl\|\,|x|^{-b}\bigl(|u|^{2\sigma}+|v|^{2\sigma}\bigr)\bigr\|_{L^{m,\tilde s}}
		\;\|u-v\|_{L^{r,2(2\sigma+1)}} .
	\end{align*}
	Applying \eqref{mstilenorm} and \Cref{lemlorentzbasic}(ii) with $2<2(2\sigma+1)$ gives \eqref{lipofnonlinearinspace}.
\end{proof}

\subsection{The space-time estimate}
We now establish the space-time estimate for the nonlinear term, which plays a crucial role in the subsequent analysis.
\begin{prop}\label{propspacetime}
	Let $(p,r)$ and $\theta$ be as in \eqref{defpr}-\eqref{deftheta} and let $T>0$. There is
	$C=C(d,b,\sigma,\lambda)>0$ such that for all $u,v\in L^{p}(0,T;L^{r,2}(\mathbb{R}^{d}))$,
	\begin{align}
	\bigl\|\mathcal{N}(u)\bigr\|_{L^{p'}(0,T;L^{r',2})}
	&\le C\,T^{\theta}\,\|u\|^{2\sigma+1}_{L^{p}(0,T;L^{r,2})},\label{nonlinearitybound}\\
	\bigl\|\mathcal{N}(u)-\mathcal{N}(v)\bigr\|_{L^{p'}(0,T;L^{r',2})}
	&\le C\,T^{\theta}\Bigl(\|u\|^{2\sigma}_{L^{p}(0,T;L^{r,2})} +\|v\|^{2\sigma}_{L^{p}(0,T;L^{r,2})}\Bigr)\|u-v\|_{L^{p}(0,T;L^{r,2})}.\label{lipnonlinearitybound}
	\end{align}
\end{prop}

\begin{proof}
	Since $\theta=1-\frac{2\sigma+2}{p}$ by \Cref{lemexponents}(iii), we have
	\begin{align}\label{timeholder}
		\frac{2\sigma+1}{p}+\theta=\frac{2\sigma+1}{p}+1-\frac{2\sigma+2}{p}=1-\frac{1}{p}=\frac{1}{p'} .
	\end{align}
	Hence, by
	\Cref{lemnonlinearspace} and the H\"older inequality in time with exponents
	$\frac{p}{2\sigma+1}$ and $\frac{1}{\theta}$,
	\begin{align*}
		\bigl\|\mathcal{N}(u)\bigr\|_{L^{p'}(0,T;L^{r',2})}
		&\le C\Bigl(\int_{0}^{T}\|u(t)\|^{(2\sigma+1)p'}_{L^{r,2}}dt\Bigr)^{1/p'}\\
		&\le C\Bigl(\int_{0}^{T}\|u(t)\|^{p}_{L^{r,2}}dt\Bigr)^{\frac{2\sigma+1}{p}}
		\Bigl(\int_{0}^{T}1\,dt\Bigr)^{\theta}
		= C\,T^{\theta}\|u\|^{2\sigma+1}_{L^{p}(0,T;L^{r,2})}.
	\end{align*}
	This proves \eqref{nonlinearitybound}. For \eqref{lipnonlinearitybound}, we first invoke
	\Cref{lemlipspace} and then apply the H\"older inequality in time with
	\begin{align*}
		\frac{1}{p'}=\frac{2\sigma}{p}+\frac{1}{p}+\ \theta.
	\end{align*}
Hence, we obtain
\begin{align*}
	\bigl\|\mathcal{N}(u)-\mathcal{N}(v)\bigr\|_{L^{p'}(0,T;L^{r',2})} &\le C \bigl\|1\bigr\|_{L^{\frac{1}{\theta}}(0,T)} \bigl\|\big(\|u\|^{2\sigma}_{L^{r,2}}+\|v\|^{2\sigma}_{L^{r,2}}\big)\bigr\|_{L^{\frac{p}{2\sigma}}(0,T)} \bigl\|\|u-v\|_{L^{r,2}}\bigr\|_{L^{p}(0,T)}\\
 &\le C\,T^{\theta}\Bigl(\|u\|^{2\sigma}_{L^{p}(0,T;L^{r,2})} +\|v\|^{2\sigma}_{L^{p}(0,T;L^{r,2})}\Bigr)\|u-v\|_{L^{p}(0,T;L^{r,2})}.
\end{align*}
This completes the proof of \eqref{lipnonlinearitybound}.
\end{proof}

\section{The truncated equation}\label{sectruncated}
In this section, we first introduce a suitable truncation of the nonlinear term and establish existence and uniqueness for the resulting equation. Next, we show that, for a sufficiently large truncation threshold, the solutions of the truncated and the original equation coincide.

Let $\vartheta\in C_{0}^{\infty}(\mathbb{R})$ with $\operatorname{supp}\vartheta\subset(-2,2)$ and
\begin{align*}
	\vartheta(x)=
	\begin{cases*}
		1, & $x\in[-1,1]$,\\
		\in[0,1], & $x\in[-2,2]$,\\
		0, & otherwise,
	\end{cases*}
\end{align*}
and for $R>0$ set $\vartheta_{R}(x):=\vartheta(x/R)$, so that $\vartheta_{R}\equiv1$ on $[-R,R]$ and
$\vartheta_{R}\equiv0$ outside $[-2R,2R]$. We consider the truncated equation
\begin{align}\label{truncated}
	\left\{
	\begin{aligned}
		du(t)&=-\bigl[\mathrm{i}Au(t)+\mathrm{i}\vartheta_{R}\bigl(\|u\|_{\XT{t}}\bigr)\mathcal{N}(u(t))
		+\beta u(t)+\mathfrak{b}(u(t))\bigr]dt\\
		&\quad-\mathrm{i}B(u(t))d\mathcal{W}_1(t)-\mathrm{i}G(u(t))d\mathcal{W}_2(t),\qquad t>0,\\
		u(0)&=u_{0}.
	\end{aligned}
	\right.
\end{align}
and the corresponding mild form,
\begin{align}\label{mildtrunc}
	u(t)&=S(t)u_{0}-\mathrm{i}\int_{0}^{t}S(t-s)\vartheta_{R}\bigl(\|u\|_{\XT{s}}\bigr)\mathcal{N}(u(s))\,ds
	-\int_{0}^{t}S(t-s)\bigl[\beta u(s)+\mathfrak{b}(u(s))\bigr]ds\nonumber\\
	&\quad-\mathrm{i}\int_{0}^{t}S(t-s)\bigl[B(u(s))d\mathcal{W}_1(s)+G(u(s))d\mathcal{W}_2(s)\bigr].
\end{align}
To establish the well-posedness of the truncated equation, we need the following theorem.
\begin{thm}\label{liptruncthm}
	Let $(p,r)$ and $\theta$ be as in \eqref{defpr}-\eqref{deftheta} and set
	\begin{align*}
		I_{R}(u)(t):=\int_{0}^{t}S(t-s)\,\vartheta_{R}\bigl(\|u\|_{\XT{s}}\bigr)\,\mathcal{N}(u(s))\,ds .
	\end{align*}
	Then for every $T>0$ the map $I_{R}$ maps $\mathcal{X}({T})$ into itself and, with
	\begin{align}\label{GammaR}
		\Gamma_{R}:=C\,(1+R)^{2\sigma},
	\end{align}
	 for all $u,u_{1},u_{2}\in\mathcal{X}({T})$, we have the following
	\begin{align}
		\|I_{R}(u)\|_{\XT{T}}&\le \Gamma_{R}\,T^{\theta}\,\|u\|_{\XT{T}},\label{boundtrunc}\\
		\|I_{R}(u_{1})-I_{R}(u_{2})\|_{\XT{T}}&\le \Gamma_{R}\,T^{\theta}\,\|u_{1}-u_{2}\|_{\XT{T}},\label{liptrunc}
	\end{align}
	where $C=C(d,b,\sigma,\lambda,\vartheta)$ does not depend on $R$, $T$, $u$, $u_{1}$ or $u_{2}$.
\end{thm}

\begin{proof}
	Fix $T>0$ and $u\in\XT{T}$, and define the stopping time
	\begin{align*}
		\tau:=\inf\bigl\{t\ge0:\ \|u\|_{\XT{t}}>2R\bigr\}\wedge T .
	\end{align*}
	Since $t\mapsto\|u\|_{\XT{t}}$ is non-decreasing and $\vartheta_{R}$ vanishes outside $[-2R,2R]$,
	we have $\vartheta_{R}(\|u\|_{\XT{t}})=0$ for $t\ge\tau$.

	By \Cref{Strichartz}(ii), applied once with the admissible pair
	$(\infty,2)$ and once with $(p,r)$, and with $(\tilde p,\tilde\varrho)=(p,r)$ in both cases,
	\begin{align}\label{boundftrunc}
		\|I_{R}(u)\|_{\XT{T}}
		\le C\bigl\|\vartheta_{R}\bigl(\|u\|_{\XT{\cdot}}\bigr)\mathcal{N}(u)\bigr\|_{L^{p'}(0,T;L^{r',2})}
		\le C\bigl\|\mathcal{N}(u)\bigr\|_{L^{p'}(0,\tau;L^{r',2})},
	\end{align}
	using $0\le\vartheta_{R}\le1$ and the vanishing of $\vartheta_{R}$ past $\tau$. By
	\eqref{nonlinearitybound} and $\|u\|_{\XT{\tau}}\le2R$,
	\begin{align*}
		\bigl\|\mathcal{N}(u)\bigr\|_{L^{p'}(0,\tau;L^{r',2})}
		\le CT^{\theta}\|u\|^{2\sigma+1}_{L^{p}(0,\tau;L^{r,2})}
		\le CT^{\theta}\|u\|^{2\sigma}_{\XT{\tau}}\,\|u\|_{\XT{T}}
		\le C(2R)^{2\sigma}T^{\theta}\|u\|_{\XT{T}}.
	\end{align*}
	This completes the proof of equation \eqref{boundtrunc}.

 Let $u_{1},u_{2}\in\XT{T}$ and define the stopping times $$\tau_{i}:=\inf\{t\ge0:\|u_{i}\|_{\XT{t}}>2R\}\wedge T,\quad i=1,2.$$ 
 Without loss of generality, we consider $\tau_{1}\le\tau_{2}$, so $[0,T]=[0,\tau_{1}]\cup[\tau_{1},\tau_{2}]\cup[\tau_{2},T]$. For simplicity, we write $F_{i}:=\vartheta_{R}(\|u_{i}\|_{\XT{\cdot}})$. We now consider the three time intervals separately. On $[0,\tau_{1}]$, we split
 \begin{align*}
 	F_{1}\mathcal{N}(u_{1})-F_{2}\mathcal{N}(u_{2})=(F_{1}-F_{2})\mathcal{N}(u_{1})+F_{2}\bigl(\mathcal{N}(u_{1})-\mathcal{N}(u_{2})\bigr)
 \end{align*}
  On $[\tau_{1},\tau_{2}]$, we have $F_{1}=0$. This implies
  \begin{align*}
  	F_{1}\mathcal{N}(u_{1})-F_{2}\mathcal{N}(u_{2})=-F_{2}\mathcal{N}(u_{2}).
  \end{align*}
	On the third subinterval $[\tau_{2},T]$, both $F_{1}$ and $F_{2}$ vanish, so the integrand is identically zero there.
	Hence, as in \eqref{boundftrunc}, \Cref{Strichartz}(ii) yields
	\begin{align}\label{l}
		\|I_{R}(u_{1})-I_{R}(u_{2})\|_{\XT{T}}
		&\le C\bigl\|\vartheta_{R}\bigl(\|u_{1}\|_{\XT{\cdot}}\bigr)\mathcal{N}(u_{1})
		-\vartheta_{R}\bigl(\|u_{2}\|_{\XT{\cdot}}\bigr)\mathcal{N}(u_{2})\bigr\|_{L^{p'}(0,T;L^{r',2})}\nonumber\\
		&\le \underbrace{C\bigl\|\bigl\{\vartheta_{R}(\|u_{1}\|_{\XT{\cdot}})
			-\vartheta_{R}(\|u_{2}\|_{\XT{\cdot}})\bigr\}\mathcal{N}(u_{1})\bigr\|_{L^{p'}(0,\tau_{1};L^{r',2})}}_{=:I^{R}_{1}}\nonumber\\
		&\quad+\underbrace{C\bigl\|\vartheta_{R}(\|u_{2}\|_{\XT{\cdot}})
			\bigl\{\mathcal{N}(u_{1})-\mathcal{N}(u_{2})\bigr\}\bigr\|_{L^{p'}(0,\tau_{1};L^{r',2})}}_{=:I^{R}_{2}}\nonumber\\
		&\quad+\underbrace{C\bigl\|\vartheta_{R}(\|u_{2}\|_{\XT{\cdot}})\mathcal{N}(u_{2})
			\bigr\|_{L^{p'}(\tau_{1},\tau_{2};L^{r',2})}}_{=:I^{R}_{3}} .
	\end{align}
	We estimate the three terms in turn, using the following lemma.

	\begin{lem}\label{lemI1}
		Let $u_{1},u_{2}\in\XT{T}$ and $w\in L^{p'}(0,T;L^{r',2}(\mathbb{R}^{d}))$. Then
		\begin{align*}
	\bigl\|\bigl\{\vartheta_{R}(\|u_{1}\|_{\XT{\cdot}})-\vartheta_{R}(\|u_{2}\|_{\XT{\cdot}})\bigr\}w
	\bigr\|_{L^{p'}(0,T;L^{r',2})}\le c_{R}\,\|u_{1}-u_{2}\|_{\XT{T}}\,\|w\|_{L^{p'}(0,T;L^{r',2})},
		\end{align*}
		with $c_{R}:=\|\vartheta_{R}'\|_{L^{\infty}(\mathbb{R})}=R^{-1}\|\vartheta'\|_{L^{\infty}(\mathbb{R})}$.
	\end{lem}
	\begin{proof}
		Applying the mean value theorem to $\vartheta_{R}$ and using the reverse triangle inequality, we achieve
		\begin{align*}
			\bigl|\vartheta_{R}(\|u_{1}\|_{\XT{t}})-\vartheta_{R}(\|u_{2}\|_{\XT{t}})\bigr|
			\le\|\vartheta_{R}'\|_{L^{\infty}}\bigl|\,\|u_{1}\|_{\XT{t}}-\|u_{2}\|_{\XT{t}}\bigr|
			\le c_{R}\|u_{1}-u_{2}\|_{\XT{T}},
		\end{align*}
		for every $t\in[0,T]$. Raising to the power $p'$, integrating in $t$ and taking the $p'$-th root
		gives the claim.
		Therefore,
		\begin{align*}
	\bigl\|\bigl\{\vartheta_{R}(\|u_{1}\|_{\XT{\cdot}})-\vartheta_{R}(\|u_{2}\|_{\XT{\cdot}})\bigr\}w
	\bigr\|_{L^{p'}(0,T;L^{r',2})}&\le \bigl\|c_{R}\|u_{1}-u_{2}\|_{\XT{T}}w\bigr\|_{L^{p'}(0,T;L^{r',2})}\\
	& \le c_{R}\,\|u_{1}-u_{2}\|_{\XT{T}}\,\|w\|_{L^{p'}(0,T;L^{r',2})}.
		\end{align*}
		This completes the proof of \Cref{lemI1}.
	\end{proof}

	\noindent \emph{{Term $I^{R}_{1}$}:} Applying \Cref{lemI1} on $[0,\tau_{1}]$ with $w=\mathcal{N}(u_{1})$, and using
	\eqref{nonlinearitybound} together with $\|u_{1}\|_{\XT{\tau_{1}}}\le 2R$, we obtain
	\begin{align}\label{estimateir1}
		I^{R}_{1}\le C\,c_{R}\,T^{\theta}\|u_{1}\|^{2\sigma+1}_{\XT{\tau_{1}}}\,\|u_{1}-u_{2}\|_{\XT{T}}
		&\le C\,R^{-1}(2R)^{2\sigma+1}\,T^{\theta}\,\|u_{1}-u_{2}\|_{\XT{T}}\nonumber\\
		&\le \tfrac{1}{3}\Gamma_{R}\,T^{\theta}\,\|u_{1}-u_{2}\|_{\XT{T}}.
	\end{align}

	\noindent \emph{{Term $I^{R}_{2}$}:} Application of $0\le\vartheta_{R}\le1$, \eqref{lipnonlinearitybound} and
	$\|u_{i}\|_{\XT{\tau_{1}}}\le2R$ gives
	\begin{align}\label{estimateir2}
		I^{R}_{2}
		&\le C\bigl\|\mathcal{N}(u_{1})-\mathcal{N}(u_{2})\bigr\|_{L^{p'}(0,\tau_{1};L^{r',2})}\nonumber\\
		&\le C\,T^{\theta}\Bigl(\|u_{1}\|^{2\sigma}_{\XT{\tau_{1}}}+\|u_{2}\|^{2\sigma}_{\XT{\tau_{1}}}\Bigr)
		\|u_{1}-u_{2}\|_{\XT{T}}\nonumber\\
		&\le 2C(2R)^{2\sigma}\,T^{\theta}\,\|u_{1}-u_{2}\|_{\XT{T}}
		\le \tfrac{1}{3}\Gamma_{R}\,T^{\theta}\,\|u_{1}-u_{2}\|_{\XT{T}} .
	\end{align}

	\noindent \emph{{Term $I^{R}_{3}$}:} On $(\tau_{1},\tau_{2})$, we have $\vartheta_{R}(\|u_{1}\|_{\XT{t}})=0$. Hence, on this interval,
	$$\vartheta_{R}(\|u_{2}\|_{\XT{\cdot}})=\vartheta_{R}(\|u_{2}\|_{\XT{\cdot}})-\vartheta_{R}(\|u_{1}\|_{\XT{\cdot}}).$$
	We apply \Cref{lemI1} as in the estimate for $I^{R}_{1}$.
	Using $\|u_{2}\|_{\XT{\tau_{2}}}\le2R$, we obtain
	\begin{align}\label{estimateir3}
		I^{R}_{3}\le C\,c_{R}\,T^{\theta}\|u_{2}\|^{2\sigma+1}_{\XT{\tau_{2}}}\,\|u_{1}-u_{2}\|_{\XT{T}}
		\le\tfrac{1}{3}\Gamma_{R}\,T^{\theta}\,\|u_{1}-u_{2}\|_{\XT{T}}.
	\end{align}
	Combining \eqref{estimateir1}-\eqref{estimateir3} with \eqref{l} gives 
	\begin{align*}
			\|I_{R}(u_{1})-I_{R}(u_{2})\|_{\XT{T}} &\le \Gamma_{R}\,T^{\theta}\,\|u_{1}-u_{2}\|_{\XT{T}}.
	\end{align*}
		This completes the proof of equation \eqref{liptrunc}
\end{proof}

\section{Well-posedness of the truncated equation}\label{secwellposedtrunc}
Now, we establish the well-posedness of the truncated equation \eqref{truncated}.
\begin{thm}\label{wellposetrunc}
	Let $q\ge2$, $T>0$, $u_{0}\in L^{\infty}(\Omega,\mathcal{F}_{0};H)$, and let $b,\sigma$ satisfy
	\eqref{brange}. Under the assumptions of \Cref{Ass}, for every $R>0$, there exists a unique
	\begin{align*}
		u^{R}\in\mathbb{X}:=L^{q}\bigl(\Omega;\mathcal{X}({T})\bigr)
	\end{align*}
	satisfying \eqref{mildtrunc}, $\mathbb{P}$-a.s., for all $t\in[0,T]$.
\end{thm}

The proof occupies the rest of this section. To begin with, we obtain a solution on a small deterministic time
interval $[0,T_{0}]$ by the Banach fixed point theorem, and then extend it to $[0,T]$ by induction.
Write
\begin{align*}
	\mathbb{X}_{1}:=L^{q}\bigl(\Omega;L^{\infty}(0,T;H)\bigr),\qquad
	\mathbb{X}_{2}:=L^{q}\bigl(\Omega;L^{p}(0,T;L^{r,2}(\mathbb{R}^{d}))\bigr),
\end{align*}
so that $\mathbb{X}=\mathbb{X}_{1}\cap\mathbb{X}_{2}$. More generally, for $0\le a<c\le T$ we write
$\mathbb{X}_{[a,c]}:=L^{q}\bigl(\Omega;\mathcal{X}({[a,c]})\bigr)$, and for convenience, we use the notation
$\mathbb{X}_{c}:=\mathbb{X}_{[0,c]}$, in particular, $\mathbb{X}=\mathbb{X}_{T}$.

\subsection{Estimates for the linear and stochastic terms}\label{subseclinear}
In this subsection, we derive estimates for the stochastic integrals appearing in the mild formulation \eqref{mildtrunc}, which will be used in establishing the local well-posedness of the truncated equation. First, let us set
\begin{align}\label{kbg}
	K_{B}(u)(t):=\int_{0}^{t}S(t-s)B(u(s))\,d\mathcal{W}_1(s),\qquad
	K_{G}(u)(t):=\int_{0}^{t}S(t-s)G(u(s))\,d\mathcal{W}_2(s).
\end{align}

\begin{lem}\label{lemKBKG}
	Write $\|B\|:=\|B\|_{\mathscr{L}(H,\mathcal{L}_{2}(Y_1,H))}$. For $u\in\mathbb{X}$ and $T\geq 0$,
	\begin{align}
		\|K_{B}(u)\|_{\mathbb{X}_{1}}+\|K_{B}(u)\|_{\mathbb{X}_{2}}
		&\le C\,T^{1/2}\,\|B\|\,\|u\|_{\mathbb{X}_{1}},\label{kbx1}\\
		\|K_{G}(u)\|_{\mathbb{X}_{1}}+\|K_{G}(u)\|_{\mathbb{X}_{2}}
		&\le C\,T^{1/2}\,\|\tilde G\|_{\infty}\,\|u\|_{\mathbb{X}_{1}}.\label{kgx1}
	\end{align}
Moreover
	\begin{align}
		\|K_{B}(u_{1})-K_{B}(u_{2})\|_{\mathbb{X}}
		&\le C\,T^{1/2}\,\|B\|\,\|u_{1}-u_{2}\|_{\mathbb{X}_{1}},\label{dif12}\\
		\|K_{G}(u_{1})-K_{G}(u_{2})\|_{\mathbb{X}}
		&\le C\,T^{1/2}\,L_{G}\,\|u_{1}-u_{2}\|_{\mathbb{X}_{1}}.\label{dif34}
	\end{align}
\end{lem}

\begin{proof}
For the $\mathbb{X}_{1}$-bounds, define $M_{B}(t):=\int_{0}^{t}S(-s)B(u(s))\,d\mathcal{W}_{1}(s)$. This is a continuous $H$-valued martingale, and the group property gives $K_{B}(u)(t)=S(t)M_{B}(t)$. By unitarity, $\|K_{B}(u)(t)\|_{H}=\|M_{B}(t)\|_{H}$. The Hilbert-space Burkholder-Davis-Gundy inequality (\Cref{BDG}) therefore applies to $M_{B}$, rather than directly to the convolution $K_{B}(u)$. For instance,
	\begin{align*}
		\|K_{B}(u)\|^{q}_{\mathbb{X}_{1}}
		&=\mathbb{E}\sup_{t\in[0,T]}\Bigl\|\int_{0}^{t}S(t-s)B(u(s))d\mathcal{W}_1(s)\Bigr\|^{q}_{H}
		\le C\,\mathbb{E}\Bigl(\int_{0}^{T}\|B(u(s))\|^{2}_{\mathcal{L}_{2}(Y_1,H)}ds\Bigr)^{q/2}\\
		&\le C\|B\|^{q}\,T^{q/2}\|u\|^{q}_{\mathbb{X}_{1}},
	\end{align*}
	and similarly for $K_{G}$, using the linear growth \eqref{LgG}, whose constant term vanishes.

	The bounds in $\mathbb{X}_{2}$ are \Cref{stochStrichartz} applied with $\Phi=B(u)$ (resp. $\Phi=G(u)$),
	together with $$\|B(u)\|_{L^{2}(0,T;\mathcal{L}_{2}(Y_1,H))}\le T^{1/2}\|B\|\,\|u\|_{L^{\infty}(0,T;H)},$$ and
	$$\|G(u)\|_{L^{2}(0,T;\mathcal{L}_{2}(Y_2,H))}\le T^{1/2}\|\tilde G\|_{\infty}\|u\|_{L^{\infty}(0,T;H)}.$$
	The constant $C(p,r,q)$ in \eqref{stochstrich} does not depend on $T$, so the whole factor $T^{1/2}$
	is supplied by these two inequalities. The difference estimates follow as $B$ is linear and $G$ is
	Lipschitz by \eqref{LipG}. The final assertion follows by running the same argument on $[a,c]$ and, for a
	stopping time, by \Cref{stochStrichartz}(ii), i.e.\ \eqref{stochstrichstop}, for the
	$\mathbb{X}_{2}$-part, and by optional stopping in \Cref{BDG} for the $\mathbb{X}_{1}$-part.
\end{proof}

Finally, for the Stratonovich correction $\mathfrak{b}$ of \eqref{defb} we have, by
\eqref{BmH} and \Cref{Strichartz}(ii) with $(\tilde p,\tilde\varrho)=(\infty,2)$,
\begin{align}\label{bx}
	\Bigl\|\int_{0}^{\cdot}S(\cdot-s)\mathfrak{b}(u(s))\,ds\Bigr\|_{\mathbb{X}}
	\le C\bigl\|\mathfrak{b}(u)\bigr\|_{L^{q}(\Omega;L^{1}(0,T;H))}
	\le C\,T\,\|u\|_{\mathbb{X}_{1}}.
\end{align}
Similarly, for the damping term $\beta u$, we have
\begin{align}\label{betaterm}
		\Bigl\|\int_{0}^{\cdot}S(\cdot-s)\beta u(s)\,ds\Bigr\|_{\mathbb{X}}
	\le C\beta T\,\|u\|_{\mathbb{X}_{1}}.
\end{align}

\subsection{Local existence}\label{subseclocal}
Define, for $u\in\mathbb{X}$,
\begin{align}\label{Iomildtrunc}
	\mathfrak{L}(u,u_{0})(t)&:=S(t)u_{0}
	-\mathrm{i}\int_{0}^{t}S(t-s)\vartheta_{R}\bigl(\|u\|_{\XT{s}}\bigr)\mathcal{N}(u(s))\,ds
	-\int_{0}^{t}S(t-s)\bigl[\beta u(s)+\mathfrak{b}(u(s))\bigr]ds\nonumber\\
	&\quad-\mathrm{i}\int_{0}^{t}S(t-s)\bigl[B(u(s))d\mathcal{W}_1(s)+G(u(s))d\mathcal{W}_2(s)\bigr] .
\end{align}
Collecting \Cref{liptruncthm}, \Cref{lemKBKG}, \eqref{bx} and \eqref{betaterm}, and using
\Cref{Strichartz}(i) for the initial term, we obtain 
\begin{align}\label{Ioxb}
	\|\mathfrak{L}(u,u_{0})\|_{\mathbb{X}}
	&\le C\|u_{0}\|_{L^{q}(\Omega;H)}
	+\Bigl[\Gamma_{R}T^{\theta}+C(1+\beta)T
	+CT^{1/2}\bigl(\|B\|+\|\tilde G\|_{\infty}\bigr)\Bigr]\|u\|_{\mathbb{X}},
\end{align}
and
\begin{align}\label{Iodifx}
	\|\mathfrak{L}(u_{1},u_{0})-\mathfrak{L}(u_{2},u_{0})\|_{\mathbb{X}}
	&\le\Bigl[\Gamma_{R}T^{\theta}+C(1+\beta)T
	+CT^{1/2}\bigl(\|B\|+L_{G}\bigr)\Bigr]\|u_{1}-u_{2}\|_{\mathbb{X}} .
\end{align}
Both right-hand sides depend linearly on $\|u\|_{\mathbb{X}}$. We now set
	\begin{align}\label{T0}
	T_{0}&=T_{0}\bigl(R,\beta,\|B\|,\|\tilde G\|_{\infty},L_{G}\bigr)\nonumber\\
	&:=\min\Bigl\{ \bigl(6\Gamma_{R}\bigr)^{-1/\theta},\ \bigl(6C(1+\beta)\bigr)^{-1},\
	\bigl(6C(1+\|B\|+\|\tilde G\|_{\infty}+L_{G})\bigr)^{-2}\Bigr\},
\end{align}
so that each of the quantities above is at most $\frac{1}{2}$ on $[0,T_{0}]$. Choose
$M:=6C\bigl(1+\|u_{0}\|_{L^{\infty}(\Omega;H)}\bigr)$ and consider the closed ball
\begin{align*}
	V^{\mathbb{X}}_{M,T_{0}}:=\bigl\{u\in\mathbb{X}:\ \|u\|_{\mathbb{X}}\le M\bigr\}
\end{align*}
over the interval $[0,T_{0}]$. Since $\|u_{0}\|_{L^{q}(\Omega;H)}\le\|u_{0}\|_{L^{\infty}(\Omega;H)}$, we
deduce from \eqref{Ioxb} that $\mathfrak{L}(\cdot,u_{0})$ maps $V^{\mathbb{X}}_{M,T_{0}}$ into itself, and
from \eqref{Iodifx} that
\begin{align}\label{contractionmathfracI}
	\|\mathfrak{L}(u_{1},u_{0})-\mathfrak{L}(u_{2},u_{0})\|_{\mathbb{X}}
	\le\tfrac{1}{2}\|u_{1}-u_{2}\|_{\mathbb{X}},\qquad u_{1},u_{2}\in V^{\mathbb{X}}_{M,T_{0}} .
\end{align}
Therefore, $\mathfrak{L}(\cdot,u_{0})$ is an $\frac{1}{2}$-contraction on $V^{\mathbb{X}}_{M,T_{0}}$, where $T_0$ depends on $R,\sigma, \beta$ but independent of the initial data $u_0$. We need the hypothesis $u_{0}\in L^{\infty}(\Omega;H)$ only to make the radius $M$ of the ball finite.
By the Banach fixed point theorem, $\mathfrak{L}(\cdot,u_{0})$ has a unique fixed point
$\overline u\in V^{\mathbb{X}}_{M,T_{0}}$, which is the unique mild solution of \eqref{truncated} on
$[0,T_{0}]$.

\subsection{Global existence for the truncated equation}\label{subsecglobaltrunc}
In this subsection, we extend the local existence of solution to global existence. Let $j:=\lfloor T/T_{0}\rfloor+1$. Suppose that, for some $k\in\{1,\dots,j\}$, there is
$u^{R}_{k}\in\mathbb{X}_{kT_{0}}$, as defined in \Cref{subseclinear}, with $u^{R}_{k}=\mathfrak{L}(u^{R}_{k},u_{0})$ on $[0,kT_{0}]$. Define
the shifted cut-off $\Theta^{R}_{k}(u)(t):=\vartheta_{R}(\phi(u)(t))$ with
\begin{align}\label{definedphi}
	\phi(u)(t)&:=\Bigl(\|u^{R}_{k}\|^{p}_{L^{p}(0,kT_{0};L^{r,2})}+\|u\|^{p}_{L^{p}(0,t;L^{r,2})}\Bigr)^{1/p}
	+\max\Bigl\{\sup_{0\le s\le kT_{0}}\|u^{R}_{k}(s)\|_{H},\ \sup_{0\le s\le t}\|u(s)\|_{H}\Bigr\}.
\end{align}

\begin{prop}\label{philip}
	For $u_{1},u_{2}\in\XT{T}$ and all $t\in[0,T]$,
	$|\phi(u_{1})(t)-\phi(u_{2})(t)|\le C\|u_{1}-u_{2}\|_{\XT{T}}$.
\end{prop}

\begin{proof}
We consider $\phi=\phi_{1}+\phi_{2}$ with
	\begin{align*}
		&\phi_1(u)(t):=\Bigl(\|u^{R}_{k}\|^{p}_{L^{p}(0,kT_{0};L^{r,2})}+\|u\|^{p}_{L^{p}(0,t;L^{r,2})}\Bigr)^{1/p}\\
		&\phi_2(u)(t):=\max\Bigl\{\sup_{0\le s\le kT_{0}}\|u^{R}_{k}(s)\|_{H},\ \sup_{0\le s\le t}\|u(s)\|_{H}\Bigr\}.
	\end{align*}
 Every norm on a normed space is Lipschitz continuous and both $\phi_1(u)$ and $\phi_2(u)$ are compositions of norms with the maps $u\mapsto u|_{[0,t]}$. Therefore, using the definition \eqref{normassum}, the elementary inequalities $|(a^{p}+x^{p})^{1/p}-(a^{p}+y^{p})^{1/p}|\le|x-y|$ and $|\max\{a,x\}-\max\{a,y\}|\le|x-y|$, we have $\phi_{1}$ is Lipschitz on $L^{p}(0,T;L^{r,2})$ and $\phi_{2}$ is Lipschitz on $C([0,T];H)$. 
\end{proof}

\noindent For $s\ge0$, we set $\widehat{\mathcal{W}}^{k}_{i}(s):=\mathcal{W}_{i}(kT_{0}+s)-\mathcal{W}_{i}(kT_{0})$,
$i=1,2$. These are, respectively, a cylindrical $Y_{1}$-valued and a standard real Wiener process with
respect to the shifted filtration $\widehat{\mathcal{F}}^{k}_{s}:=\mathcal{F}_{kT_{0}+s}$. Define
$\mathfrak{L}_{k}$ on $[0,T_{0}]$ by
\begin{align*}
	\mathfrak{L}_{k}(u)(t)&:=S(t)u^{R}_{k}(kT_{0})
	-\mathrm{i}\int_{0}^{t}S(t-s)\Theta^{R}_{k}(u)(s)\mathcal{N}(u(s))\,ds
	-\int_{0}^{t}S(t-s)\bigl[\beta u(s)+\mathfrak{b}(u(s))\bigr]ds\nonumber\\
	&\quad-\mathrm{i}\int_{0}^{t}S(t-s)\bigl[B(u(s))\,d\widehat{\mathcal{W}}^{k}_{1}(s)
	+G(u(s))\,d\widehat{\mathcal{W}}^{k}_{2}(s)\bigr].
\end{align*}

\begin{prop}\label{lipjkr}
	Let $(p,r)$ and $\theta$ be as in \eqref{defpr}-\eqref{deftheta} and set
	\begin{align*}
		\mathfrak{L}^{R}_{k}(u)(t):=\mathrm{i}\int_{0}^{t}S(t-s)\Theta^{R}_{k}(u)(s)\mathcal{N}(u(s))\,ds.
	\end{align*}
Then, for every $T>0$, $\mathfrak{L}^{R}_{k}$ maps $\mathcal{X}({T})$ into itself and for $u_1, u_2\in \mathcal{X}({T})$, we have
	\begin{align*}
		\|\mathfrak{L}^{R}_{k}(u)\|_{\XT{T}}\le\Gamma_{R}\,T^{\theta}\,\|u\|_{\XT{T}},
		\qquad
		\|\mathfrak{L}^{R}_{k}(u_{1})-\mathfrak{L}^{R}_{k}(u_{2})\|_{\XT{T}}
		\le \Gamma_{R}\,T^{\theta}\,\|u_{1}-u_{2}\|_{\XT{T}},
	\end{align*}
	with $\Gamma_{R}$ as in \eqref{GammaR}
\end{prop}

\begin{proof}
	The proof is immediate from \Cref{liptruncthm} and \Cref{philip}. The proof of \Cref{liptruncthm} uses only two properties of $\vartheta_{R}$. That it is non-decreasing in $t$ and
	Lipschitz in $u$ with respect to $\|\cdot\|_{\XT{T}}$. By \Cref{philip}, $\phi$ also has both properties.
\end{proof}
\begin{proof}[Proof of \Cref{wellposetrunc}:]
Repeating the argument of \Cref{subseclocal} with \Cref{lipjkr} as a substitute in place of \Cref{liptruncthm} shows that $\mathfrak{L}_{k}$ is a $\tfrac12$-contraction on $\mathbb{X}_{[kT_{0},(k+1)T_{0}]}$. Let $\overline u^{R}_{k+1}$ be its unique fixed point and set 
\begin{align*}
	u^{R}_{k+1}(t):=
	\begin{cases}
		u^{R}_{k}(t), & t\in[0,kT_{0}],\\
		\overline u^{R}_{k+1}(t-kT_{0}), & t\in[kT_{0},(k+1)T_{0}].
	\end{cases}
\end{align*}
Using $\vartheta_{R}(\|u^{R}_{k}\|_{\XT{s}})=\vartheta_{R}(\|u^{R}_{k+1}\|_{\XT{s}})$ for $s\in[0,kT_{0}]$
and $\Theta^{R}_{k}(\overline u^{R}_{k+1})(s)=\vartheta_{R}(\|u^{R}_{k+1}\|_{\XT{kT_{0}+s}})$ for
$s\in[0,T_{0}]$, together with the semigroup property $S(\hat t)S(kT_{0})=S(kT_{0}+\hat t)$ and the change
of variable $s\mapsto kT_{0}+s$, one can easily verify that $u^{R}_{k+1}$ satisfies
\eqref{mildtrunc} on $[0,(k+1)T_{0}]$. After $j$ steps, we obtain the unique solution $u^{R}:=u^{R}_{j}$ on $[0,T]$ to the truncated equation \eqref{truncated}. This completes the proof of \Cref{wellposetrunc}. 
\end{proof}
\section{ The mass identity and removal of the truncation}\label{secremoval}
In this section, first we establish an exact mass identity for the truncated equation by means of a suitable Yosida approximation, which yields a pathwise bound on the $L^2$-norm independent of the truncation parameter. Then, we combine this bound with a subcritical bootstrap to obtain a uniform space-time estimate. Next, we show that, for sufficiently large truncation thresholds, the truncated solution coincides with the solution of the original equation. This allows us to remove the truncation and complete the proof of \Cref{mainthm}.
\subsection{An exact mass identity}
In this subsection, we first establish an exact mass identity for the truncated solution by applying the It\^o formula to a suitable Yosida regularization of the $L^2$-mass functional. Consequently, we obtain a pathwise bound for the $L^2$-mass that is independent of the truncation parameter.
\begin{lem}\label{lemmassapprox}
Let $v\in C([0,T];H)\cap L^{p}(0,T;L^{r,2})$ satisfy the truncated mild equation.
For $J_{\epsilon}=(I+\epsilon A)^{-1}$, the process $v_{\epsilon}=J_{\epsilon}v$ is an
$H$-valued strong semimartingale and It\^o's formula applies to $\|v_{\epsilon}\|_{H}^{2}$.
Moreover, after integration in time, every term converges in probability, uniformly on compact
time intervals, to the corresponding duality pairing in the formal $L^{2}$ identity.
\end{lem}
\begin{proof}
The resolvent is self-adjoint on $H$, commutes with $A$ and $S(t)$. It is uniformly bounded on
$L^{s,2}$ for $1<s<\infty$, and converges strongly on these spaces.We first show that, for each fixed $\epsilon>0$ the resolvent maps $L^{r',2}$ into $H$. Indeed, we can write $J_{\epsilon}f=G_{\epsilon}*f$ with $G_{\epsilon}(x)=\epsilon^{-d/2}G_{1}(x/\sqrt\epsilon)$, where $G_{1}$ is the Bessel kernel of order
$2$. Hence $G_{1}(x)=O(|x|^{2-d})$ as $|x|\to0$ (logarithmically for $d=2$, boundedly for $d=1$).
$G_{1}$ also decays exponentially at infinity, so $G_{\epsilon}\in L^{m}(\mathbb R^{d})$ whenever
$1/m>1-2/d$. Using O'Neil's convolution inequality \cite{ONeil1963}, we obtain
\begin{align}\label{Jeps}
 \|J_{\epsilon}f\|_{H}=\|G_{\epsilon}*f\|_{L^{2,2}}
 \le C\,\|G_{\epsilon}\|_{L^{m}}\,\|f\|_{L^{r',2}},
 \qquad \frac1m=1+\frac12-\frac1{r'},
\end{align} 
Since $1<r'<2$, we have $1<m<\infty$. The required integrability $1/m>1-2/d$ is equivalent to $\frac32-\frac1{r'}>1-\frac2d$, which can also be written as $d\bigl(\frac1{r'}-\frac12\bigr)<2$. This holds since admissibility gives $d\bigl(\frac1{r'}-\frac12\bigr)=2/p<1$. We do not use the sharp Lorentz-Sobolev embedding here. This embedding would give $L^{\varrho^{*},2}$ with $1/\varrho^{*}=1/r'-2/d<\tfrac12$, and no embedding $L^{\varrho^{*},2}\hookrightarrow L^{2}$ is available on $\mathbb R^{d}$. Thus $J_{\epsilon}\mathcal{N}(v)\in L^{p'}(0,T;H)$. Since $p'>1$ and $T<\infty$, this regularized nonlinear drift also belongs to $L^{1}(0,T;H)$, while $AJ_{\epsilon}v$, $J_{\epsilon}v$, $J_{\epsilon}\mathfrak b(v)$ and the regularised noise coefficients have the integrability required by the infinite dimensional It\^o formula \cite[Th. 4.32]{DaPrato+Zabzyk_2014}. For the nonlinear drift, self-adjointness of $J_{\epsilon}$ gives
\begin{align*}
	 \int_{0}^{t}(J_{\epsilon}\mathcal{N}(v),J_{\epsilon}v)_{H}\,ds
	=\int_{0}^{t}\langle\mathcal{N}(v),J_{\epsilon}^{2}v\rangle_{L^{r',2},L^{r,2}}\,ds,
\end{align*}
which converges by H\"older in time, the nonlinear estimate \eqref{nonlinearitybound}, and strong
convergence of $J_{\epsilon}^{2}v$ in $L^{p}_{t}L^{r,2}_{x}$. The linear terms converge by strong
$H$-convergence and \eqref{BmH}. The quadratic-variation sums converge by dominated convergence theorem using $\sum_m\|B_m\|^{2}<\infty$ and \eqref{LgG}. The Burkholder-Davis-Gundy inequality (\Cref{BDG}) gives that the martingale terms converge uniformly in time in probability. Finally
$\|J_{\epsilon}v\|_{H}^{2}\to\|v\|_{H}^{2}$ uniformly in time. This follows from the compactness of
$v([0,T])\subset H$ together with the uniform boundedness of $J_{\epsilon}$.
\end{proof}

\begin{prop}[Mass identity]\label{massidentity}
	Let $u^{R}$ be the solution of \eqref{truncated} given by \Cref{wellposetrunc}. Then, $\mathbb{P}$-a.s.,
	the function $t\mapsto\|u^{R}(t)\|^{2}_{H}$ is absolutely continuous on $[0,T]$ and solves the
	\emph{deterministic} differential equation
	\begin{align}\label{massODE}
		\frac{d}{dt}\,y^{R}(t)=\Bigl(\tilde G\bigl(y^{R}(t)\bigr)^{2}-2\beta\Bigr)\,y^{R}(t)
		\quad\text{for a.e.\ }t\in[0,T],
		\qquad y^{R}(0)=\|u_{0}\|^{2}_{H},
	\end{align}
	where $y^{R}(t):=\|u^{R}(t)\|^{2}_{H}$. Consequently
	\begin{align}\label{massbound}
		\sup_{t\in[0,T]}\|u^{R}(t)\|^{2}_{H}
		\le\|u_{0}\|^{2}_{H}\,e^{\left(\|\tilde G\|^{2}_{\infty}-2\beta\right)_{+}T}
		\le\|u_{0}\|^{2}_{L^{\infty}(\Omega;H)}\,e^{\left(\|\tilde G\|^{2}_{\infty}-2\beta\right)_{+}T}
		=:\mathcal{M},
		\qquad\mathbb{P}\text{-a.s.},
	\end{align}
	and $\mathcal{M}$ is independent of $R$ and of $\omega$. If $2\beta\ge\|\tilde G\|^{2}_{\infty}$ then
	$y^{R}$ is non-increasing and $\|u^{R}(t)\|_{H}\le\|u_{0}\|_{H}$ for all $t$.
\end{prop}

\begin{proof}
Let us define $u^{\epsilon}:=J_{\epsilon}u^{R}$. By \Cref{lemmassapprox}, $u^{\epsilon}$ solves the
	truncated equation in the strong sense. therefore, we can apply the infinite dimensional It\^o formula (\cite[Th. 4.32]{DaPrato+Zabzyk_2014}) to $\|u^{\epsilon}\|^{2}_{H}$. This gives, for every $t\in[0,T]$,
	\begin{align}\label{regito}
		\|u^{\epsilon}(t)\|^{2}_{H}
		&=\|J_{\epsilon}u_{0}\|^{2}_{H}
		+2\int_{0}^{t}\bigl(-\mathrm{i}Au^{\epsilon},u^{\epsilon}\bigr)_{H}ds
		+2\int_{0}^{t}\vartheta_{R}\bigl(-\mathrm{i}J_{\epsilon}\mathcal{N}(u^{R}),u^{\epsilon}\bigr)_{H}ds
		\nonumber\\
		&\quad-2\beta\int_{0}^{t}\|u^{\epsilon}\|^{2}_{H}ds
		-2\int_{0}^{t}\bigl(J_{\epsilon}\mathfrak{b}(u^{R}),u^{\epsilon}\bigr)_{H}ds\nonumber\\
		&\quad+\int_{0}^{t}\Bigl[\sum_{m=1}^{\infty}\bigl\|J_{\epsilon}B_{m}u^{R}\bigr\|^{2}_{H}
		+\bigl\|J_{\epsilon}G(u^{R})\bigr\|^{2}_{\mathcal{L}_{2}(Y_2,H)}\Bigr]ds\nonumber\\
		&\quad+2\sum_{m=1}^{\infty}\int_{0}^{t}\bigl(-\mathrm{i}J_{\epsilon}B_{m}u^{R},u^{\epsilon}\bigr)_{H}d\beta^{1}_{m}
		+2\int_{0}^{t}\bigl(-\mathrm{i}J_{\epsilon}G(u^{R}),u^{\epsilon}\bigr)_{H}\,d\beta^{2} .
	\end{align}
	The second term on the right-hand side vanishes identically. Indeed $J_{\epsilon}$ maps $H$ onto
	$D(A)$, so $u^{\epsilon}(s)\in D(A)$, and $A$ is self-adjoint for the Hermitian inner product, whence
	$\langle Au^{\epsilon},u^{\epsilon}\rangle\in\mathbb{R}$ and
	$(-\mathrm{i}Au^{\epsilon},u^{\epsilon})_{H}
	=\operatorname{Re}(-\mathrm{i}\langle Au^{\epsilon},u^{\epsilon}\rangle)=0$.

	Letting $\epsilon\downarrow0$ in \eqref{regito} and using \Cref{lemmassapprox} term by term, we
	obtain, $\mathbb{P}$-a.s.\ and for every $t\in[0,T]$,
	\begin{align}\label{limitito}
		\|u^{R}(t)\|^{2}_{H}
		&=\|u_{0}\|^{2}_{H}
		+2\int_{0}^{t}\vartheta_{R}\bigl(-\mathrm{i}\mathcal{N}(u^{R}),u^{R}\bigr)_{H}ds
		-2\beta\int_{0}^{t}\|u^{R}\|^{2}_{H}ds
		-2\int_{0}^{t}\bigl(\mathfrak{b}(u^{R}),u^{R}\bigr)_{H}ds\nonumber\\
		&\quad+\int_{0}^{t}\Bigl[\sum_{m=1}^{\infty}\bigl\|{-\mathrm{i}}B_{m}u^{R}\bigr\|^{2}_{H}
		+\bigl\|G(u^{R})\bigr\|^{2}_{\mathcal{L}_{2}(Y_2,H)}\Bigr]ds\nonumber\\
		&\quad+2\sum_{m=1}^{\infty}\int_{0}^{t}\bigl(-\mathrm{i}B_{m}u^{R},u^{R}\bigr)_{H}d\beta^{1}_{m}
		+2\int_{0}^{t}\bigl(-\mathrm{i}G(u^{R}),u^{R}\bigr)_{H}\,d\beta^{2}.
	\end{align}
      For the term involving $\mathcal{N}(u^{R})$, we use the duality pairing between $L^{r',2}$ and $L^{r,2}$ given in \Cref{lemlorentzbasic}(v). We treat the remaining terms one by one.

	\noindent First, we consider the nonlinear term. Since $\lambda\in\mathbb{R}$ and $|x|^{-b}|u|^{2\sigma}$ is real-valued,
	\begin{align*}
		\bigl(-\mathrm{i}\mathcal{N}(u),u\bigr)_{H}
		=\operatorname{Re}\Bigl(-\mathrm{i}\lambda\int_{\mathbb{R}^{d}}|x|^{-b}|u|^{2\sigma}u\,\overline u\,dx\Bigr)
		=\lambda\operatorname{Re}\Bigl(-\mathrm{i}\int_{\mathbb{R}^{d}}|x|^{-b}|u|^{2\sigma+2}dx\Bigr)=0 ,
	\end{align*}
	the integral being a real number. The cut-off factor $\vartheta_{R}$ is a scalar and does not affect
	this. The integral is finite by the \emph{endpoint} Lorentz duality estimate
	\Cref{lemlorentzbasic}(v). Indeed, \eqref{keyidentity} yields $\frac{b}{d}+\frac{2\sigma+2}{r}=1$, by \Cref{lemlorentzbasic}(ii)-(iii) and \Cref{lemweight},
	$|u|^{2\sigma+2}\in L^{d/(d-b),1}$ and
	\begin{align*}
	\int_{\mathbb R^d}|x|^{-b}|u|^{2\sigma+2}\,dx
	&\le \big\||x|^{-b}\big\|_{L^{d/b,\infty}}
	\big\||u|^{2\sigma+2}\big\|_{L^{d/(d-b),1}}\le C\|u\|_{L^{r,2}}^{2\sigma+2}<\infty,
	\end{align*}
	 for $t$ a.e. 
	 \noindent Now, we deal the Stratonovich correction term. By \eqref{Bmselfadjoint},
	$$(\mathfrak{b}(u),u)_{H}=\frac12\sum_{m}(B_{m}^{2}u,u)_{H}=\frac12\sum_{m}\|B_{m}u\|^{2}_{H},$$ while
	$\sum_{m}\|-\mathrm{i}B_{m}u\|^{2}_{H}=\sum_{m}\|B_{m}u\|^{2}_{H}$. Hence
	$-2(\mathfrak{b}(u),u)_{H}+\sum_{m}\|-\mathrm{i}B_{m}u\|^{2}_{H}=0$. 
	
	\noindent For the martingale terms, \eqref{Bmselfadjoint} implies $\langle B_{m}u,u\rangle\in\mathbb{R}$.
	Hence
	\begin{align*}
		(-\mathrm{i}B_{m}u,u)_{H}=\operatorname{Re}(-\mathrm{i}\langle B_{m}u,u\rangle)=0.
	\end{align*}
	Furthermore, using \eqref{perticularg} and \eqref{L2R}, we obtain
	\begin{align*}
		\bigl(-\mathrm{i}G(u),u\bigr)_{H}
		=\tilde G\bigl(\|u\|^{2}_{H}\bigr)\operatorname{Re}\bigl(-\mathrm{i}\langle u,u\rangle\bigr)
		=\tilde G\bigl(\|u\|^{2}_{H}\bigr)\operatorname{Re}\bigl(-\mathrm{i}\|u\|^{2}_{H}\bigr)=0 .
	\end{align*}

\noindent Again by \eqref{L2R},
	$\|G(u)\|^{2}_{\mathcal{L}_{2}(Y_{2},H)}=\tilde G(\|u\|^{2}_{H})^{2}\|u\|^{2}_{H}$. Thus, \eqref{limitito} gives \eqref{massODE}, since all the stochastic terms vanish. In particular $y^{R}=\|u^{R}\|^{2}_{H}$ is, pathwise, absolutely continuous and satisfies an ordinary differential equation. Since $\tilde G$ is bounded, $\bigl(y^{R}\bigr)'\le(\|\tilde G\|^{2}_{\infty}-2\beta)y^{R}$ a.e., and Gronwall's lemma gives \eqref{massbound}. Moreover, if $2\beta\ge\|\tilde G\|^{2}_{\infty}$, then right-hand side is less than or equal to $0$, so $y^{R}$ is non-increasing. This completes the proof of \Cref{massidentity}.
\end{proof}

\subsection{A uniform Strichartz bound}

The mass bound \eqref{massbound} controls only $\|u^{R}\|_{L^{\infty}(0,T;H)}$. To remove the truncation we
must control the norm $\|u^{R}\|_{\XT{T}}$, by a bound that remains uniform in $R$.

First, the cut-off has no impact on $[0,\tau_R]$. The nonlinear bound \eqref{nonlinearitybound} has power $2\sigma+1$. If we take the $L^q(\Omega)$ norm directly, the right-hand side requires an
$L^{q(2\sigma+1)}(\Omega)$ norm of the path norm. In \Cref{sectruncated}, the cut-off gives a linear bound. Here we introduce a second family of stopping times at a level $K$ chosen independently of $R$. Up to these stopping times, we use $K$ to bound the factor with power $2\sigma$. This leaves a linear bound in the path norm.

Second, we use \eqref{massbound} in each subinterval. At every initial time $t_j$, we have
$\|u^{R}(t_j)\|_H\le \mathcal{M}^{1/2}$, where $\mathcal{M}$ is independent of $R$ and $\omega$. Therefore, we can choose the length $\delta$ of each subinterval using only $\mathcal{M}$ and $K$. Here, subcriticality is used. Since $\theta>0$, we can make $\delta^\theta$ as small as needed by choosing $\delta$ small.

\begin{prop}\label{propuniformbound}
	Let $q\ge2$, let $u^{R}$ be the solution of \eqref{truncated} given by \Cref{wellposetrunc}, let
	$\mathcal{M}$ be as in \eqref{massbound}, and set
	\begin{align}\label{tauR}
		\tau_{R}:=\inf\bigl\{t\ge0:\ \|u^{R}\|_{\XT{t}}>R\bigr\}\wedge T .
	\end{align}
	There are constants $\mathcal{K}\ge1$ and $c_{0}\ge1$, depending only on
	$T,\mathcal{M},\beta,\|B\|,\|\tilde G\|_{\infty},d,b,\sigma$ and $q$, and \emph{not} on $R$, with the
	following property. Let $K\ge\mathcal{K}$, put
	\begin{align}\label{deltachoice}
		\delta:=\min\Bigl\{1,\ \bigl(8C_{\ast}K^{2\sigma}\bigr)^{-1/\theta},\ \delta_{1}\Bigr\},\qquad
		N:=\lceil T/\delta\rceil,\qquad t_{j}:=j\delta,
	\end{align}
	where $C_{\ast}$ is the constant of \eqref{nonlinearitybound} and $\delta_{1}\in(0,1]$ is fixed in
	\eqref{delta1} below, and define the stopping times
	\begin{align}\label{varrhoj}
		\varrho_{j}:=\inf\bigl\{t\ge t_{j}:\ \|u^{R}\|_{\XT{[t_{j},t]}}>K\bigr\}\wedge t_{j+1}\wedge\tau_{R},
		\qquad j=0,\dots,N-1 .
	\end{align}
	Then, writing $\mathcal{A}_{K}:=\bigcap_{j=0}^{N-1}\bigl\{\varrho_{j}=t_{j+1}\wedge\tau_{R}\bigr\}$,
	\begin{enumerate}
		\item[(i)] $\bigl\|u^{R}\bigr\|_{L^{q}\left(\Omega;\,\mathcal{X}({[t_{j},\varrho_{j}]})\right)}\le\mathcal{K}$
		for every $j=0,\dots,N-1$;
		\item[(ii)] $\displaystyle\mathbb{P}\bigl(\mathcal{A}^{c}_{K}\bigr)\le N\,\mathcal{K}^{q}K^{-q}
		\le c_{0}\,\mathcal{K}^{q}\,K^{\frac{2\sigma}{\theta}-q}$;
		\item[(iii)] on $\mathcal{A}_{K}$, 
		$\bigl\|u^{R}\bigr\|_{\XT{\tau_{R}}}\le\Lambda(K):=c_{0}\,K^{\,a}$, where
		$a:=1+\dfrac{2\sigma}{p\theta}$.
	\end{enumerate}
\end{prop}

\begin{proof}
	On $[0,\tau_{R}]$ we have $\|u^{R}\|_{\XT{t}}\le R$, hence $\vartheta_{R}(\|u^{R}\|_{\XT{t}})=1$, so
	$u^{R}$ satisfies the untruncated mild equation \eqref{milditoeq} there. Restarting it at $t_{j}$ and
	using the semigroup property, we obtain for $t\in[t_{j},\varrho_{j}]$
	\begin{align*}
		u^{R}(t)&=S(t-t_{j})u^{R}(t_{j})
		-\int_{t_{j}}^{t}S(t-s)\bigl[\mathrm{i}\mathcal{N}(u^{R}(s))+\beta u^{R}(s)+\mathfrak{b}(u^{R}(s))\bigr]ds\\
		&\quad-\mathrm{i}\int_{t_{j}}^{t}S(t-s)\bigl[B(u^{R}(s))\,d\mathcal{W}_1(s)
		+G(u^{R}(s))\,d\beta^{2}(s)\bigr].
	\end{align*}
	Applying \Cref{Strichartz} on $[t_{j},t]$ pathwise, together with \eqref{nonlinearitybound} and
	\eqref{bx}, gives
	\begin{align}\label{piecewise}
		\|u^{R}\|_{\XT{[t_{j},t]}}
		\le C\|u^{R}(t_{j})\|_{H}
		+C_{\ast}\delta^{\theta}\|u^{R}\|^{2\sigma+1}_{\XT{[t_{j},t]}}
		+C(1+\beta)\delta\,\|u^{R}\|_{\XT{[t_{j},t]}}
		+\mathcal{Z}_{j}(t),
	\end{align}
	where $\mathcal{Z}_{j}(t)$ denotes the $\XT{[t_{j},t]}$-norm of the two stochastic convolutions
	$K_{B}(u^{R})$ and $K_{G}(u^{R})$ started at $t_{j}$.

   By \eqref{varrhoj} and the continuity of $t\mapsto\|u^{R}\|_{\XT{[t_{j},t]}}$, we
	have $\|u^{R}\|_{\XT{[t_{j},t]}}\le K$ for every $t\le\varrho_{j}$. Therefore, using the second entry of \eqref{deltachoice}
	\begin{align*}
		C_{\ast}\delta^{\theta}\|u^{R}\|^{2\sigma+1}_{\XT{[t_{j},t]}}
		\le C_{\ast}\delta^{\theta}K^{2\sigma}\,\|u^{R}\|_{\XT{[t_{j},t]}}
		\le\dfrac18\,\|u^{R}\|_{\XT{[t_{j},t]}},\qquad t\le\varrho_{j}.
	\end{align*}
	 On $[t_{j},\varrho_{j}]$ the nonlinear term is linear in $\|u^{R}\|$, so that taking
	$L^{q}(\Omega)$-norms below is justified. Such a $\delta$ exists precisely because $\theta>0$, that is, $\sigma<\frac{2-b}{d}$. Fixing in addition
	\begin{align}\label{delta1}
		\delta_{1}:=\min\Bigl\{1,\ \bigl(8C(1+\beta)\bigr)^{-1}\Bigr\}
	\end{align}
	and absorbing the two linear terms into the left-hand side of \eqref{piecewise}, we obtain the pathwise
	bound
	\begin{align}\label{pathwisebound}
		\|u^{R}\|_{\XT{[t_{j},\varrho_{j}]}}\le 2C\|u^{R}(t_{j})\|_{H}+2\,\mathcal{Z}_{j}(\varrho_{j})
		\le 2C\mathcal{M}^{1/2}+2\,\mathcal{Z}_{j}(\varrho_{j}),
	\end{align}
	where \eqref{massbound} was used in the last step.

	\emph{(i).} Take $L^{q}(\Omega)$-norms in \eqref{pathwisebound} and apply \Cref{lemKBKG} on
	$[t_{j},\varrho_{j}]$. More precisely, extend
	$\mathbf1_{[t_j,\varrho_j]}B(u^{R})$ and
	$\mathbf1_{[t_j,\varrho_j]}G(u^{R})$ by zero to $[t_j,t_{j+1}]$ and apply \eqref{stochstrichstop} of \Cref{stochStrichartz}(ii).  Progressive measurability is preserved and the constant
	is independent of $j$ and $\varrho_j$. Using \eqref{massbound} once more,
	\begin{align*}
		\bigl\|\mathcal{Z}_{j}(\varrho_{j})\bigr\|_{L^{q}(\Omega)}
		\le C\delta^{1/2}\bigl(\|B\|+\|\tilde G\|_{\infty}\bigr)
		\bigl\|u^{R}\bigr\|_{L^{q}(\Omega;L^{\infty}(t_{j},\varrho_{j};H))}
		\le C\bigl(\|B\|+\|\tilde G\|_{\infty}\bigr)\mathcal{M}^{1/2}.
	\end{align*}
	Hence (i) holds with
	\begin{align}\label{defK}
		\mathcal{K}:=1+2C\mathcal{M}^{1/2}+2C\bigl(\|B\|+\|\tilde G\|_{\infty}\bigr)\mathcal{M}^{1/2},
	\end{align}
	which is independent of $R$, $K$, $j$, $\delta$.

	\emph{(ii).} If $\varrho_{j}<t_{j+1}\wedge\tau_{R}$ then, by continuity of
	$t\mapsto\|u^{R}\|_{\XT{[t_{j},t]}}$, necessarily $\|u^{R}\|_{\XT{[t_{j},\varrho_{j}]}}=K$. Therefore, Chebyshev's
	inequality and (i) gives
	$$\mathbb{P}\bigl(\varrho_{j}<t_{j+1}\wedge\tau_{R}\bigr)\le K^{-q}\,\mathbb{E}
	\|u^{R}\|^{q}_{\XT{[t_{j},\varrho_{j}]}}\le\mathcal{K}^{q}K^{-q},$$ and a union bound over the $N$ values
	of $j$ yields the first inequality. For the second, \eqref{deltachoice} gives, for $K\ge1$,
	\begin{align}\label{Nbound}
		N\le 1+\frac{T}{\delta}\le 1+\frac{T}{\delta_{1}}+T\,(8C_{\ast})^{1/\theta}K^{\frac{2\sigma}{\theta}}
		\le c_{0}\,K^{\frac{2\sigma}{\theta}} .
	\end{align}

	\emph{(iii).} On $\mathcal{A}_{K}$ we have $\varrho_{j}=t_{j+1}\wedge\tau_{R}$ and, by the definition
	\eqref{varrhoj}, $\|u^{R}\|_{\XT{I_{j}}}\le K$ for every $j$, where
	$I_{j}:=[t_{j},t_{j+1}]\cap[0,\tau_{R}]$. Since
	$\|u^{R}\|_{L^{\infty}(0,\tau_{R};H)}=\max_{j}\|u^{R}\|_{L^{\infty}(I_{j};H)}\le K$ and
	$\|u^{R}\|^{p}_{L^{p}(0,\tau_{R};L^{r,2})}=\sum_{j}\|u^{R}\|^{p}_{L^{p}(I_{j};L^{r,2})}\le NK^{p}$, we get
	\begin{align*}
		\|u^{R}\|_{\XT{\tau_{R}}}\le\bigl(1+N^{1/p}\bigr)K
		\le\Bigl(1+c_{0}^{1/p}K^{\frac{2\sigma}{p\theta}}\Bigr)K
		\le c_{0}\,K^{1+\frac{2\sigma}{p\theta}}
	\end{align*}
	by \eqref{Nbound}, after enlarging $c_{0}$ if necessary.
\end{proof}

\subsection{Proof of \texorpdfstring{\Cref{mainthm}}{Theorem 1.1}}
 In this section, we establish the proof of the main result \Cref{mainthm}. For that, we first need to establish the following lemma:
\begin{lem}\label{lemlocalunique}
	Let $u,v\in L^{q}(\Omega;\mathcal{X}({T}))$ be mild solutions of \eqref{milditoeq} with the same initial
	datum. Let $\varrho\le T$ be a stopping time for which there is a deterministic $n<\infty$ with
	\begin{align}\label{localbound}
		\|u\|_{\XT{t}}+\|v\|_{\XT{t}}\le n\qquad\text{for all }t\le\varrho,\ \mathbb{P}\text{-a.s.}
	\end{align}
Then $u=v$ on $[0,\varrho]$, $\mathbb{P}$-a.s. In particular, the conclusion holds for
	$$\zeta_{n}:=\inf\{t\ge0:\|u\|_{\XT{t}}+\|v\|_{\XT{t}}>n\}\wedge T.$$
\end{lem}

\begin{proof}
	Since $\theta>0$, we may choose a deterministic $h\in(0,1]$ with
	\begin{align}\label{hchoice}
		2C\,n^{2\sigma}h^{\theta}+C(1+\beta)h+C\,h^{1/2}\bigl(\|B\|+L_{G}\bigr)\le\tfrac12.
	\end{align}
We set $t_{k}:=kh$. Suppose, inductively, that $u=v$ on $[0,t_{k}\wedge\varrho]$. 

This is clearly true for $k=0$. We now consider both mild equations starting from $t_k\wedge\varrho$, with their equal initial values. As explained in \Cref{subsecglobaltrunc}, the shifted Wiener processes are Brownian motions with respect to the shifted filtration. Moreover, by \eqref{localbound}, we have $\|u\|_{\XT{t}},\|v\|_{\XT{t}}\le n$ for all
$t\le\varrho$, a.s. Hence,	\eqref{lipnonlinearitybound} gives
	\begin{align*}
		\bigl\|\mathcal{N}(u)-\mathcal{N}(v)\bigr\|_{L^{p'}\left(t_{k}\wedge\varrho,\ t_{k+1}\wedge\varrho;\,L^{r',2}\right)}
		\le 2C\,n^{2\sigma}h^{\theta}\,
		\|u-v\|_{\XT{[t_{k}\wedge\varrho,\ t_{k+1}\wedge\varrho]}},
	\end{align*}
an estimate that is linear in $\|u-v\|$. Combining it with \Cref{Strichartz}, \eqref{bx} and the
	stopped-interval form of \Cref{lemKBKG}, and using \eqref{hchoice}, we obtain
	\begin{align*}
		X:=\bigl\|u-v\bigr\|_{L^{q}\left(\Omega;\,\XT{[t_{k}\wedge\varrho,\ t_{k+1}\wedge\varrho]}\right)}
		\le\frac12\,X .
	\end{align*}
	Since $u,v\in L^{q}(\Omega;\mathcal{X}({T}))$, we have $X<\infty$, and therefore $X=0$. Hence $u=v$ on
	$[0,t_{k+1}\wedge\varrho]$. After $\lceil T/h\rceil$ steps, we have $u=v$ on $[0,\varrho]$.
\end{proof}
Finally, we are ready to prove the main result, Theorem \ref{mainthm}. We do so in the following steps. 
\begin{proof}[Proof of \Cref{mainthm}]
	Fix $q\ge2$ and set $a:=1+\dfrac{2\sigma}{p\theta}$ as in \Cref{propuniformbound}. We choose an auxiliary exponent
	\begin{align}\label{barq}
	\bar q:=\dfrac{2\sigma}{\theta}+aq+1 .
	\end{align}
	Since $a>1$, we have $\bar q>q\ge2$ and $\bar q>\frac{2\sigma}{\theta}$. We apply \Cref{wellposetrunc} and \Cref{propuniformbound} with $\bar q$ in place of $q$. Let $\mathcal{K},c_{0}$ be the corresponding constants and write $\Lambda(K)=c_{0}K^{a}$. Since $\mathbb{P}$ is a probability measure, we have
	\begin{align*}
		L^{\bar q}(\Omega;X)\hookrightarrow L^{q}(\Omega;X).
	\end{align*}
	Thus, it is enough to construct the solution with moment exponent $\bar q$. Choose $R$ large enough so that
	\begin{align*}
		K_{R}:=\left(\frac{R}{2c_{0}}\right)^{1/a}\ge\mathcal{K}.
	\end{align*}
	Apply \Cref{propuniformbound} with $K=K_{R}$. By \Cref{propuniformbound}(iii), on the event $\mathcal{A}_{K_{R}}$, we have 
	$$\|u^{R}\|_{\XT{\tau_{R}}}\le\Lambda(K_{R})=R/2<R.$$ If $\tau_{R}<T$, continuity of $t\mapsto\|u^{R}\|_{\XT{t}}$ and \eqref{tauR} imply $\|u^{R}\|_{\XT{\tau_{R}}}=R$. Hence $\tau_{R}=T$ on $\mathcal{A}_{K_{R}}$. Since $K_{R}\to\infty$ and $\bar q>\frac{2\sigma}{\theta}$ by \eqref{barq}, from \Cref{propuniformbound}(ii), we have
	\begin{align}\label{tauRprob}
		\mathbb{P}\bigl(\tau_{R}<T\bigr)\le\mathbb{P}\bigl(\mathcal{A}^{c}_{K_{R}}\bigr)
		\le c_{0}\mathcal{K}^{\bar q}\,K_{R}^{\frac{2\sigma}{\theta}-\bar q}
		\xrightarrow[R\to\infty]{}0.
	\end{align}

\noindent Let $R_{1}<R_{2}$ and put
	\begin{align*}
		\varrho:=\inf\bigl\{t\ge0:\ \|u^{R_{1}}\|_{\XT{t}}\vee\|u^{R_{2}}\|_{\XT{t}}>R_{1}\bigr\}\wedge T .
	\end{align*}
	On $[0,\varrho]$ both cut-offs equal $1$, so $u^{R_{1}}$ and $u^{R_{2}}$ are both mild solutions of the
	untruncated equation \eqref{milditoeq} there. Moreover $\|u^{R_{1}}\|_{\XT{t}}\vee\|u^{R_{2}}\|_{\XT{t}}\le R_{1}$
	for $t\le\varrho$, so \Cref{lemlocalunique} applies with $n=2R_{1}$.
	It follows that
	\begin{align*}
		u^{R_{1}}=u^{R_{2}}\quad\text{on }[0,\varrho].
	\end{align*}
	In particular,
	\begin{align*}
		\|u^{R_{2}}\|_{\XT{t}}
		=\|u^{R_{1}}\|_{\XT{t}},
		\qquad t<\varrho.
	\end{align*}
	Hence, the two processes reach the level $R_{1}$ at the same time.
	Therefore,
	\begin{align*}
		\varrho=\tau_{R_{1}}, \qquad u^{R_{1}}=u^{R_{2}}\quad\text{on }[0,\tau_{R_{1}}],
	\quad \hbox{and} \quad
		\tau_{R_{1}}\le\tau_{R_{2}}.
	\end{align*}
	Thus, the family $(\tau_{R})_{R>0}$ is non-decreasing. By \eqref{tauRprob}, we can choose a deterministic sequence $R_{k}\uparrow\infty$ such that
	\begin{align*}
		\mathbb{P}\bigl(\tau_{R_{k}}<T\bigr)\le2^{-k}.
	\end{align*}
	Since $\sum_{k}\mathbb{P}(\tau_{R_{k}}<T)<\infty$, the Borel-Cantelli lemma gives a $\mathbb{P}$-null set $\mathcal{N}$ such that, for every $\omega\in\Omega\setminus\mathcal{N}$, we have $\tau_{R_{k}}(\omega)=T$ for all sufficiently large $k$. For such $\omega$, define
	\begin{align*}
		u(\cdot,\omega):=u^{R_{k}}(\cdot,\omega)
		\quad\text{on }[0,T],
	\end{align*}
	where $k$ is any sufficiently large index.  Moreover, $u$ is
	adapted, has continuous paths in $H$, and satisfies
	\eqref{milditoeq} on $[0,T]$.

	\noindent Let $\Lambda\ge\Lambda_{0}:=c_{0}\mathcal{K}^{a}$ and
	$K:=(\Lambda/c_{0})^{1/a}\ge\mathcal{K}$. Consequently, applying \Cref{propuniformbound} on $\mathcal{A}_{K}$, we get $\|u^{R}\|_{\XT{\tau_{R}}}\le\Lambda<R$. Hence $\tau_{R}=T$. This implies that $u=u^{R}$ on $[0,T]$ and $\|u\|_{\XT{T}}\le\Lambda$. Therefore
	$\bigl\{\|u\|_{\XT{T}}>\Lambda\bigr\}\subseteq\mathcal{A}^{c}_{K}$, and \Cref{propuniformbound}(ii) gives
	\begin{align}\label{tail}
		\mathbb{P}\bigl(\|u\|_{\XT{T}}>\Lambda\bigr)
		\le c_{0}\mathcal{K}^{\bar q}K^{\frac{2\sigma}{\theta}-\bar q}
		=c_{1}\,\Lambda^{-\frac{1}{a}\left(\bar q-\frac{2\sigma}{\theta}\right)},
		\qquad \Lambda\ge\Lambda_{0} .
	\end{align}
	By \eqref{barq} the exponent satisfies $\frac{1}{a}\bigl(\bar q-\frac{2\sigma}{\theta}\bigr)>q$, whence
	\begin{align}\label{momentbound}
		\mathbb{E}\,\|u\|^{q}_{\XT{T}}
		=q\int_{0}^{\infty}\Lambda^{q-1}\,\mathbb{P}\bigl(\|u\|_{\XT{T}}>\Lambda\bigr)\,d\Lambda
		\le\Lambda_{0}^{q}+q\,c_{1}\int_{\Lambda_{0}}^{\infty}
		\Lambda^{q-1-\frac{1}{a}\left(\bar q-\frac{2\sigma}{\theta}\right)}d\Lambda<\infty .
	\end{align}
	Thus $u\in L^{q}\bigl(\Omega;L^{\infty}(0,T;H)\cap L^{p}(0,T;L^{r,2})\bigr)$.

\noindent Continuity follows from the mild formulation \eqref{milditoeq}.
	The initial term is continuous by the strong continuity of $S(\cdot)$.
	The deterministic Duhamel terms are continuous in $H$ by \Cref{Strichartz}(ii). Furthermore, the nonlinear forcing lies in
	$L^{p'}(0,T;L^{r',2})$ and the bounded linear drift terms lie in $L^1(0,T;H)$. Each stochastic convolution has a continuous $H$-valued
	version. 
	 Hence $u\in L^{q}\bigl(\Omega;\mathcal{X}({T})\bigr)$.

\end{proof}

\begin{remark}\label{remopen}
	Two natural questions are left open. 
	\begin{enumerate}
		\item[(i)] In the mass-critical case $\sigma=\frac{2-b}{d}$, we have $\theta=0$. The proof of \Cref{propuniformbound} uses $\theta>0$ to make $\delta^\theta$ small. In the deterministic case, \cite{AnKim2021} treats the critical case by taking small initial data and working in Lorentz spaces. We expect a similar argument to work here when $\|u_0\|_{L^\infty(\Omega;H)}$
		is small.

		\item[(ii)] In the energy space $H^1(\mathbb{R}^d)$, differentiating the nonlinearity also produces a derivative of the singular weight. This reduces the range of admissible values of $b$.
	\end{enumerate}
	The auxiliary exponent $\bar q$ in \eqref{barq} arises from the localisation in \Cref{propuniformbound}. It remains open whether the moment bound \eqref{momentbound} also holds with $\bar q=q$.

\end{remark}

\noindent \textbf{Acknowledgements:}
The authors gratefully acknowledge Dr. Manil T. Mohan (Associate Professor, Department of Mathematics, IIT Roorkee) and Prof. Erika Hausenblas (Chair of Applied Mathematics, Montanuniversit\"at Leoben) for valuable suggestions, insightful discussions, and constructive comments that contributed to the development of this work. 

\noindent \textbf{Data availibility:} Data sharing not applicable to this article as no datasets were generated or analysed during the current study.\\

\noindent \textbf{Funding:} The second named author of the paper acknowledges the financial support provided by the Anusandhan National Research Foundation (ANRF) under the Early Career Research Grant (ECRG), Grant No. ANRF/ECRG/2025/002863/PMS.\\

\noindent \textbf{Author contributions:} Sandip Roy and Debopriya Mukherjee both authors contributed to the study conception, investigation, analysis and manuscript writing. All authors read and approved the final manuscript.\\

\noindent \textbf{\large Declarations}\\

\noindent\textbf{Conflict of interest:} The author has no competing interests to declare that are relevant to the content of this article.

\begin{appendix}
	\section{Preliminaries of Lorentz space}
We recall the following five standard properties of Lorentz spaces;
	see \cite{Grafakos,ONeil1963}.
	
	\begin{lem}[Basic properties]\label{lemlorentzbasic}
		Let $0<\varrho<\infty$ and $0<s_{1}\le s_{2}\le\infty$. Then
		\begin{enumerate}
			\item[(i)] $L^{\varrho,\varrho}(\mathbb{R}^{d})=L^{\varrho}(\mathbb{R}^{d})$ with equality of norms;
			\item[(ii)] $L^{\varrho,s_{1}}(\mathbb{R}^{d})\hookrightarrow L^{\varrho,s_{2}}(\mathbb{R}^{d})$,
			i.e.\ $\|f\|_{L^{\varrho,s_{2}}}\le C\|f\|_{L^{\varrho,s_{1}}}$;
			\item[(iii)] for $\alpha>0$ and $0<s\le\infty$,
			$\bigl\||f|^{\alpha}\bigr\|_{L^{\varrho,s}}=\|f\|^{\alpha}_{L^{\alpha\varrho,\alpha s}}$;
			\item[(iv)] for $1<\varrho<\infty$ and $1\le s<\infty$, $\bigl(L^{\varrho,s}\bigr)^{*}=L^{\varrho',s'}$;
			\item[(v)] \emph{(duality pairing)} for $1<\varrho<\infty$ and $1\le s\le\infty$,
			\begin{align}\label{lorentzduality}
				\int_{\mathbb{R}^{d}}|f(x)g(x)|\,dx\;\le\;\|f\|_{L^{\varrho,s}}\,\|g\|_{L^{\varrho',s'}} .
			\end{align}
		\end{enumerate}
	\end{lem}
	
	\begin{proof}
		(i) and (ii) are classical. For (iii), first suppose $s<\infty$. Then
		$\bigl(|f|^{\alpha}\bigr)^{*}=(f^{*})^{\alpha}$ and $\bigl(t^{1/(\alpha\varrho)}\bigr)^{\alpha s}=t^{s/\varrho}$, so $\bigl\||f|^{\alpha}\bigr\|^{s}_{L^{\varrho,s}}
		=\int_{0}^{\infty}t^{s/\varrho}\bigl(f^{*}(t)\bigr)^{\alpha s}\frac{dt}{t}
		=\|f\|^{\alpha s}_{L^{\alpha\varrho,\alpha s}}$.
	For $s=\infty$, the same identity follows by taking suprema in the definition of the weak Lorentz quasi-norm. (iv) is the duality theorem for Lorentz spaces, it fails for $s=\infty$. For (v), the Hardy-Littlewood rearrangement inequality gives $\int|fg|\le\int_{0}^{\infty}f^{*}(t)g^{*}(t)\,dt$, and the H\"older inequality on $(0,\infty)$ with the measure $dt/t$ applied to $f^{*}g^{*}=\bigl(t^{1/\varrho}f^{*}\bigr)\bigl(t^{1/\varrho'}g^{*}\bigr)t^{-1}$ yields \eqref{lorentzduality} with constant $1$. See \cite[\S1.4]{Grafakos}.
	\end{proof}
	
	\begin{lem}[Lorentz-H\"older inequality, \cite{ONeil1963}]\label{lemoneil}
		Let $1<\varrho_{1},\varrho_{2}<\infty$ and $0<s_{1},s_{2}\le\infty$ satisfy
		\begin{align*}
			\frac{1}{\varrho}=\frac{1}{\varrho_{1}}+\frac{1}{\varrho_{2}}<1,
			\qquad \frac{1}{s}=\frac{1}{s_{1}}+\frac{1}{s_{2}}.
		\end{align*}
		Then there is $C=C(\varrho_{1},\varrho_{2},s_{1},s_{2})$ with
		\begin{align}\label{oneil}
			\|fg\|_{L^{\varrho,s}}\le C\,\|f\|_{L^{\varrho_{1},s_{1}}}\|g\|_{L^{\varrho_{2},s_{2}}}
		\end{align}
		for all $f\in L^{\varrho_{1},s_{1}}(\mathbb{R}^{d})$, $g\in L^{\varrho_{2},s_{2}}(\mathbb{R}^{d})$.
	\end{lem}
	
	\begin{proof}
		See \cite[Theorem 3.4]{ONeil1963} or \cite[Exercise 1.4.19]{Grafakos}.
	\end{proof}

	The following elementary computation is the substitute for the Hardy-Littlewood-Sobolev inequality, the singular weight in the weak Lebesgue scale.
	\begin{lem}\label{lemweight}
		Let $d\ge1$ and $0<b<d$. Then $|x|^{-b}\in L^{\frac{d}{b},\infty}(\mathbb{R}^{d})$, and
		\begin{align}\label{weightnorm}
			\bigl\||x|^{-b}\bigr\|_{L^{\frac{d}{b},\infty}(\mathbb{R}^{d})}=\omega_{d}^{\,b/d},
			\qquad \omega_{d}:=\bigl|B(0,1)\bigr| .
		\end{align}
		In particular the norm is finite and explicit, whereas $|x|^{-b}\notin L^{\varrho}(\mathbb{R}^{d})$
		for every $\varrho\in[1,\infty]$.
	\end{lem}
	
	\begin{proof}
		Write $w(x)=|x|^{-b}$. For $\alpha>0$,
		$\{|w|>\alpha\}=\{|x|<\alpha^{-1/b}\}$, hence $d_{w}(\alpha)=\omega_{d}\,\alpha^{-d/b}$.
		Inverting, $w^{*}(t)=(t/\omega_{d})^{-b/d}$, so that $t^{b/d}w^{*}(t)=\omega_{d}^{b/d}$ for
		every $t>0$, which is \eqref{weightnorm}. The last assertion is \eqref{nolebesgue}, proved in
		\Cref{secintro}.
	\end{proof}
	\subsection{Deterministic Strichartz estimates in Lorentz spaces}\label{subsecstrichartz}
	\noindent Keel and Tao's abstract argument yields, through real rather than complex interpolation, a
	sharpening of the Strichartz estimates in which the Lebesgue space $L^{\varrho}$ is replaced by the
	Lorentz space $L^{\varrho,2}$; since $L^{\varrho,2}\hookrightarrow L^{\varrho,\varrho}=L^{\varrho}$ for
	$\varrho\ge2$, this is a genuine improvement. We record the version we need.
	
	\begin{prop}[Strichartz estimates in the Lorentz scale]\label{Strichartz}
		Let $\{S(t)\}_{t\in\mathbb{R}}$ be as in \Cref{subsecnotation} and let $(p,\varrho)$ and
		$(\tilde p,\tilde\varrho)$ be admissible pairs with $p,\tilde p>2$.
		\begin{enumerate}
			\item[(i)] For every $\phi\in H$ the map $t\mapsto S(t)\phi$ belongs to
			$L^{p}(\mathbb{R};L^{\varrho,2}(\mathbb{R}^{d}))\cap L^{\infty}(\mathbb{R};H)$ and
			\begin{align}\label{strichartz1}
				\|S(\cdot)\phi\|_{L^{p}(\mathbb{R};L^{\varrho,2}(\mathbb{R}^{d}))}\le C\|\phi\|_{H}.
			\end{align}
 \item[(ii)] Let $I\subseteq\mathbb{R}$ be an interval with $0\in\overline I$. For every $g\in L^{\tilde p'}(I;L^{\tilde\varrho',2}(\mathbb{R}^{d}))$ the map $\Phi^{S}_{g}(t):=\int_{0}^{t}S(t-s)g(s)\,ds$ belongs to $L^{p}(I;L^{\varrho,2}(\mathbb{R}^{d}))\cap C(\overline I;H)$ and, with $C$ independent of $I$,
			\begin{align}\label{strichartz2}
				\|\Phi^{S}_{g}\|_{L^{\infty}(\overline I;H)}
				+\|\Phi^{S}_{g}\|_{L^{p}(I;L^{\varrho,2}(\mathbb{R}^{d}))}
				\le C\|g\|_{L^{\tilde p'}(I;L^{\tilde\varrho',2}(\mathbb{R}^{d}))}.
			\end{align}
		\end{enumerate}
	\end{prop}
	
	\begin{proof}
Estimate \eqref{strichartz1} is the Lorentz refinement of the Keel-Tao theorem, i.e., the real interpolation argument in \cite[\S6]{KeelTao1998} yields the spatial $L^{\varrho,2}$ norm. The corresponding inhomogeneous estimate is obtained by the Lorentz duality ($ (L^{\tilde\varrho,2})^*=L^{\tilde\varrho',2}$; see \Cref{lemlorentzbasic}(iv)). Since $p>2>\tilde p'$, the standard estimate gives the $L^p_tL^{\varrho,2}_x$ part of \eqref{strichartz2}, i.e., the $L^\infty_tH$ part follows by duality with the homogeneous estimate for $(\tilde p,\tilde\varrho)$ (or directly by unitarity when $(\tilde p,\tilde\varrho)=(\infty,2)$). We cite \cite{AnKim2021}. To obtain continuity in $H$, approximate $g$ in $L^{\tilde p'}(I;L^{\tilde\varrho',2})$ by smooth functions compactly supported in space and time. Applying \eqref{strichartz2} with pair $(\infty,2)$, one can show uniform convergence in $H$ of the corresponding integrals. This gives a continuous representative for $\Phi^S_g$.
	\end{proof}
	
	\begin{remark}\label{remendpoint}
		The admissible pair produced in \Cref{lemexponents} always satisfies $p>2(\sigma+1)>2$, so the
		restriction $p>2$ in \Cref{Strichartz} does not work. The pair $(\infty,2)$ is admissible and
		$L^{2,2}=L^{2}=H$, so \eqref{strichartz2} with $(\tilde p,\tilde\varrho)=(\infty,2)$ recovers the
		usual energy estimate.
	\end{remark}
	
		\subsection{Stochastic Strichartz estimates in Lorentz spaces}
	
	\begin{lem}[Burkholder-Davis-Gundy inequality]\label{BDG}
		Let $M$ be an $H$-valued continuous $\mathbb{F}$-local martingale with $M_{0}=0$. For every
		$q\in(0,\infty)$ and every $\mathbb{F}$-stopping time $\tau$,
		\begin{align*}
			\mathbb{E}\sup_{t\le\tau}\|M_{t}\|_{H}^{q}\;\approx_{q}\;\mathbb{E}[M,M]^{q/2}_{\tau}
			=\mathbb{E}\langle M,M\rangle^{q/2}_{\tau}.
		\end{align*}
	\end{lem}
This is the standard Hilbert-space BDG inequality; see,
	for example, \cite{DaPrato+Zabzyk_2014}.
	
	\noindent Finally we shall need a Strichartz estimate for the stochastic convolution, again in the
	Lorentz scale. It rests on the following geometric property of the target space.

	\begin{lem}[$2$-smoothness of the Strichartz space]\label{lem2smooth}
	Let $2\le p<\infty$, $2\le\varrho<\infty$, $T>0$ and $X:=L^{p}(0,T;L^{\varrho,2}(\mathbb{R}^{d}))$. Then $X$ admits an equivalent lattice norm which is $2$-uniformly smooth, with equivalence constants independent of $T$. In particular $X$ has martingale type $2$ and type $2$, with constants independent of $T$.
	\end{lem}
	\begin{proof}
	We first verify $2$-convexity. Since $\varrho/2\ge1$, $L^{\varrho/2,1}$ is a Banach lattice. Thus, by \Cref{lemlorentzbasic}(iii), for every finite family $(f_j)$,
	\begin{align*}
	\left\|\left(\sum_j |f_j|^2\right)^{1/2}\right\|_{L^{\varrho,2}}^2
	&=\left\|\sum_j|f_j|^2\right\|_{L^{\varrho/2,1}}
	\le\sum_j\|f_j\|_{L^{\varrho,2}}^2.
	\end{align*}
	Minkowski's inequality in $L^{p/2}(0,T)$, using $p\ge2$, therefore shows that $X$ is $2$-convex with constant $1$. For finite concavity, the classical lower-estimate results for Lorentz spaces imply that $L^{\varrho,2}$ is $s$-concave for any $s>\varrho$ when $\varrho>2$; for $\varrho=2$ this is just the corresponding property of $L^2$ \cite{Creekmore1981,LTII}. Choose $s>\max\{p,\varrho\}$. Minkowski's inequality then gives
	\begin{align*}
	\left(\sum_j\|f_j\|_X^s\right)^{1/s}
	&\le \left\|\left(\sum_j\|f_j(t)\|_{L^{\varrho,2}}^s
	\right)^{1/s}\right\|_{L^p_t}\\
	&\le M_s\left\|\left(\sum_j|f_j|^s\right)^{1/s}\right\|_X,
	\end{align*}
	where $M_s$ depends on $s$ and $\varrho$, but not on $T$. Thus $X$ is $2$-convex and $s$-concave for some finite $s$, with both constants independent of $T$. The Banach-lattice renorming theorem \cite[Theorem~1.f.1]{LTII} yields an equivalent lattice norm on $X$ with modulus of smoothness of
	power type $2$ and uniform equivalence constants. The martingale type $2$ conclusion follows from uniform smoothness \cite{Pisier1975}, and martingale type $2$ implies type $2$.
	\end{proof}
	
	\begin{prop}[Stochastic Strichartz estimate]\label{stochStrichartz}
		Let $(p,\varrho)$ be an admissible pair with $2<p<\infty$, $2\le\varrho<\infty$, let $q\ge2$, $T>0$,
		let $Y$ be a separable Hilbert space, let $\mathcal{W}$ be a cylindrical $\mathbb{F}$-Wiener process on $Y$,
		and let $\Phi$ be an $\mathcal{L}_{2}(Y,H)$-valued progressively
		measurable process with
		$\mathbb{E}\|\Phi\|^{q}_{L^{2}(0,T;\mathcal{L}_{2}(Y,H))}<\infty$.
		\begin{enumerate}
			\item[(i)] The stochastic convolution
			$K_{\Phi}(t):=\int_{0}^{t}S(t-s)\Phi(s)\,d\mathcal{W}(s)$ satisfies
			\begin{align}\label{stochstrich}
				\bigl\|K_{\Phi}\bigr\|_{L^{q}(\Omega;L^{p}(0,T;L^{\varrho,2}(\mathbb{R}^{d})))}
				\le C(p,\varrho,q)\,
				\bigl\|\Phi\bigr\|_{L^{q}(\Omega;L^{2}(0,T;\mathcal{L}_{2}(Y,H)))},
			\end{align}
			with a constant $C(p,\varrho,q)$ independent of $T$.
			\item[(ii)] Let $0\le a<c\le T$ and let $\tau$ be a
			stopping time with $a\le\tau\le c$. Then, with the same constant as in \eqref{stochstrich},
			\begin{align}\label{stochstrichstop}
				\Bigl\|\int_{a}^{\cdot}S(\cdot-s)\Phi(s)\,d\mathcal{W}(s)\Bigr\|_{L^{q}(\Omega;L^{p}(a,\tau;L^{\varrho,2}(\mathbb{R}^{d})))}
				\le C(p,\varrho,q)\,
				\bigl\|\Phi\bigr\|_{L^{q}(\Omega;L^{2}(a,\tau;\mathcal{L}_{2}(Y,H)))},
			\end{align}
			where the norms over the random interval $(a,\tau)$ are taken pathwise.
		\end{enumerate}
	\end{prop}
	
	\noindent Since the constant in \eqref{stochstrich} does not depend on $T$, every power of $T$
	occurring below is extracted from the right-hand side of \eqref{stochstrich}.
	\begin{proof}
		\emph{(i).} Let $X:=L^{p}(0,T;L^{\varrho,2}(\mathbb{R}^{d}))$, equipped with the equivalent norm of
		\Cref{lem2smooth}. 		For $s\in[0,T]$ define 
		$$J_{s}:H\to X \quad \hbox{by} \quad
		(J_{s}h)(t):=\mathbf{1}_{(s,T)}(t)\,S(t-s)h.$$ Since $S(t-s)h=S(t)S(-s)h$, the homogeneous estimate
		\eqref{strichartz1} applied to $S(-s)h$ and unitarity give
		\begin{align}\label{Jsbound}
			\|J_{s}h\|_{X}\le\|S(\cdot)S(-s)h\|_{L^{p}(\mathbb{R};L^{\varrho,2})}\le C_{0}\|h\|_{H},
			\qquad s\in[0,T],\ h\in H,
		\end{align}
  with $C_{0}$ being the constant of \eqref{strichartz1}. Moreover $s\mapsto J_{s}h$ is continuous in $X$: $s\mapsto S(-s)h$ is continuous in $H$, \eqref{Jsbound} is uniform in $s$, and for fixed $F\in L^{p}(\mathbb{R};L^{\varrho,2})$ the map $s\mapsto\mathbf{1}_{(s,T)}F$ is continuous in $X$ by dominated convergence. The required progressive measurability of $\widetilde\Phi(s):=J_s\Phi(s)$ in the space $\gamma(Y,X)$ will follow from the uniform bound below.  The factor $\mathbf 1_{\{s<t\}}$ ensures that only times before $t$ are considered, so the process remains adapted. \\
	Since $X$ has type $2$, for $\Psi\in\mathcal{L}_{2}(Y,H)$ and an orthonormal basis $\{e_{m}\}$ of $Y$, by \eqref{Jsbound},
		\begin{align}\label{gammanorm}
			\|J_{s}\Psi\|_{\gamma(Y,X)}
			\le\tau_{2}(X)\Bigl(\sum_{m}\|J_{s}\Psi e_{m}\|^{2}_{X}\Bigr)^{1/2}
			\le\tau_{2}(X)\,C_{0}\,\|\Psi\|_{\mathcal{L}_{2}(Y,H)} .
		\end{align}
	For each fixed finite-rank $\Psi$, the continuity of $s\mapsto J_sh$ implies continuity of $s\mapsto J_s\Psi$ in $\gamma(Y,X)$. By the uniform estimate \eqref{gammanorm} and density of finite-rank operators in $\mathcal L_2(Y,H)$, this extends to every $\Psi\in\mathcal L_2(Y,H)$. Consequently the map $(s,\Psi)\mapsto J_s\Psi$ is jointly continuous into $\gamma(Y,X)$, and the progressive measurability of $\Phi$ implies that of $\widetilde\Phi$.\\ 
	Since $X$ is $2$-smooth (\Cref{lem2smooth}), the stochastic integral $\int_{0}^{T}\widetilde\Phi\,d\mathcal{W}$ exists in $X$ and, for every $q\in(0,\infty)$,
		\begin{align*}
			\mathbb{E}\Bigl\|\int_{0}^{T}\widetilde\Phi(s)\,d\mathcal{W}(s)\Bigr\|^{q}_{X}
			\le C_{q}\,\mathbb{E}\Bigl(\int_{0}^{T}\|\widetilde\Phi(s)\|^{2}_{\gamma(Y,X)}ds\Bigr)^{q/2}
			\le C_{q}\bigl(\tau_{2}(X)C_{0}\bigr)^{q}\,
			\mathbb{E}\|\Phi\|^{q}_{L^{2}(0,T;\mathcal{L}_{2}(Y,H))},
		\end{align*}
		by \cite{Brzezniak1997,Seidler2010} and \eqref{gammanorm}. In particular, after taking $q$th roots,
		\begin{align*}
		\left\|\int_0^T\widetilde\Phi(s)\,d\mathcal W(s)
		\right\|_{L^q(\Omega;X)}
		\le C_{p,\varrho,q}
		\|\Phi\|_{L^q(\Omega;L^2(0,T;\mathcal L_2(Y,H)))},
		\end{align*}
		where the constant is independent of $T$ by \Cref{lem2smooth}.\\
Let $g\in C_c^\infty((0,T)\times\mathbb R^d)$ and set $\ell_g(f):=\operatorname{Re}\int_0^T\!\int_{\mathbb R^d} f(t,x)\overline{g(t,x)}\,dx\,dt$. As $g$ is smooth and compactly supported in \emph{both} variables, $\ell_g$ is bounded on $X$ (by Lorentz duality) and on $L^2(0,T;H)$. The stochastic Fubini theorem applying to the unitarity of $S(t)$,
	\begin{align*}
	&\int_0^T\!\int_0^T\mathbf 1_{\{s<t\}} \big\|\langle S(t-s)\Phi(s)\,\cdot\,,g(t)\rangle
	\big\|_{\mathcal L_2(Y,\mathbb R)}^2\,dt\,ds\\
	&\qquad\le \|g\|_{L^2(0,T;H)}^2 \|\Phi\|_{L^2(0,T;\mathcal L_2(Y,H))}^2,
	\end{align*}
	and the right-hand side belongs to $L^{q/2}(\Omega)$. Bounded linear functionals commute with stochastic integration, so stochastic Fubini \cite{DaPrato+Zabzyk_2014} yields
		\begin{align*}
	\ell_{g}\Bigl(\int_{0}^{T}\widetilde\Phi\,d\mathcal{W}\Bigr)=\int_{0}^{T}\Bigl(\int_{s}^{T}\operatorname{Re}\bigl\langle S(t-s)\Phi(s)\,\cdot\,,g(t)\bigr\rangle dt\Bigr)d\mathcal{W}(s) =\int_{0}^{T}\operatorname{Re}\bigl\langle K_{\Phi}(t),g(t)\bigr\rangle dt =\ell_{g}(K_{\Phi}),
		\end{align*}
$\mathbb{P}$-a.s. We choose a countable separating family of such smooth test functions. The corresponding identities hold simultaneously outside a single null set and identify the $H$-valued convolution with the $X$-valued integral as a space-time distribution. Hence $K_{\Phi}=\int_{0}^{T}\widetilde\Phi\,d\mathcal{W}$ in $X$, $\mathbb{P}$-a.s. This proves \eqref{stochstrich} for every $q\in(0,\infty)$, in particular for every $q\ge2$. None of $C_{0}$, $\tau_{2}(X)$, $C_{q}$ depends on $T$.\\
\emph{(ii).} We put $M(t):=\int_{0}^{t}S(-s)\Phi(s)\,d\mathcal{W}(s)$, a continuous $H$-valued martingale. Then $K_{\Phi}(t)=S(t)M(t)$. The process $\Phi_{\tau}:=\mathbf{1}_{(a,\tau]}\Phi$ is progressively measurable, since $\mathbf{1}_{(a,\tau]}$ is a predictable time indicator, and $K_{\Phi_{\tau}}(t)=S(t)\bigl(M(t\wedge\tau)-M(t\wedge a)\bigr)$ for all $t$, a.s. Hence, for$t\in[a,\tau]$, $$K_{\Phi_{\tau}}(t)=S(t)\bigl(M(t)-M(a)\bigr)=\int_{a}^{t}S(t-s)\Phi(s)\,d\mathcal{W}(s).$$
Consequently, the left-hand side of \eqref{stochstrichstop} is at most $\|K_{\Phi_{\tau}}\|_{L^{q}(\Omega;L^{p}(0,c;L^{\varrho,2}))}$, which by (i) on $[0,c]$ is at most $C(p,\varrho,q)\|\Phi_{\tau}\|_{L^{q}(\Omega;L^{2}(0,c;\mathcal{L}_{2}(Y,H)))}$, and the latter equals the right-hand side of \eqref{stochstrichstop}. This completes the proof.
	\end{proof}

	\end{appendix}

\bibliography{Inhomogeneous_final}
\bibliographystyle{abbrv}

\end{document}